\documentclass[11pt]{article}

\usepackage{arxiv}
\usepackage[utf8]{inputenc}
\usepackage[T1]{fontenc}
\usepackage[pagebackref=true]{hyperref}
\usepackage{url}
\usepackage{booktabs}
\usepackage{amsfonts}
\usepackage{nicefrac}
\usepackage{microtype}
\usepackage{graphicx}
\usepackage{natbib}
\usepackage{doi}
\usepackage{dsfont}
\usepackage{amssymb,amsmath,amsfonts,amsthm}
\usepackage{verbatim}	
\usepackage{fancyvrb}
\usepackage{bbm}
\usepackage{paralist}
\usepackage{accents}
\usepackage{verbatim}
\usepackage[ruled]{algorithm2e}
\usepackage{cleveref}
\usepackage{color}
\usepackage{caption}
\usepackage{multirow}
\usepackage{stix}
\usepackage{enumitem}
\usepackage{mathtools}
\usepackage{float}
\usepackage{multicol}
\usepackage{comment}
\usepackage{import}
\usepackage{subfig}
\usepackage{tikz-cd}
\usetikzlibrary{babel}
\usepackage{accents}
\usepackage{wrapfig}
\usepackage{extarrows}

\DeclareMathOperator*{\esssup}{ess\,sup}
\DeclareMathOperator{\Vrt}{Vert}

\def\R{\mathbb{R}}
\def\N{\mathbb{N}}
\def\Z{\mathbb{Z}}
\def\Q{\mathbb{Q}}
\newcommand{\E}[1]{\mathbb{E}\left[{#1}\right]}
\newcommand{\Et}[1]{\tilde{\mathbb{E}}\left[{#1}\right]}
\newcommand{\V}[1]{\mathbb{V}\left({#1}\right)}
\newcommand{\Cov}[2]{\mathrm{Cov}\left({#1},{#2}\right)}
\def\P{\mathbb{P}}

\def\A{\mathcal{A}}
\def\W{\mathcal{W}}

\def\MM{\mathbb{M}}

\def\AA{\mathbb{A}}

\newcommand{\1}{\mathds{1}}
\newcommand{\0}{\mathbf{0}}

\newcommand{\norm}[1]{\Vert #1 \Vert}

\newcommand{\ND}{\mathds{N}(\R^d\times\AA)}

\renewcommand{\d}[1]{ \mathrm{d}{#1} }

\newcommand{\fWn}{f\big(\Delta_{W_n}\big)}
\newcommand{\fD}[1]{f\big(\Delta_{{#1}}\big)}
\renewcommand{\|}{\mid}

\newcommand{\chic}{\chi_{\mathrm{C}}^r}
\newcommand{\chivr}{\chi_{\mathrm{VR}}^r}

\newtheorem{Th}{Theorem}[section]
\newtheorem{Prop}[Th]{Proposition}
\newtheorem{Lem}[Th]{Lemma}
\newtheorem{Kor}[Th]{Corollary}
\newtheorem{Bem}[Th]{Remark}
\newtheorem{Def}[Th]{Definition}

\hypersetup{pdftoolbar=true,
	pdfmenubar=true,
	pdfpagemode=UseOutlines,
	bookmarksnumbered=true,
	linktocpage=true,
	colorlinks=true,
	citecolor=blue,
	linkcolor=blue,
	urlcolor=blue
}

\title{Betti numbers in the Random Connection Model for higher-dimensional simplicial complexes and the Boolean model}

\author{ Dominik~Pabst \\
    Institute of Stochastics\\
    Karlsruhe Institute of Technology\\
    Institute of Theoretical Physics\\
    Friedrich Alexander University Erlangen-Nuremberg
}

\renewcommand{\headeright}{}
\renewcommand{\undertitle}{}
\renewcommand{\shorttitle}{Betti numbers in the RCM for higher-dimensional simplicial complexes}

\hypersetup{
pdftitle={Betti numbers in the Random Connection Model for higher-dimensional simplicial complexes and the Boolean model},
pdfauthor={Dominik~Pabst},
pdfkeywords={Random Connection Model, Betti numbers, Central limit theorems, Boolean model, Simplicial complexes},
}

\begin{document}

\maketitle

\begingroup
\renewcommand{\thefootnote}{}
\footnotetext{© 2026. This manuscript version is made available under the CC-BY-NC-ND 4.0 license \url{https://creativecommons.org/licenses/by-nc-nd/4.0/}.}
\addtocounter{footnote}{-1}
\endgroup

\begin{abstract}
Random simplicial complexes, as generalizations of random graphs, have become increasingly popular in the literature in recent years.
In this paper, we consider a new model for a random simplicial complex that was introduced in \cite{Pabst.Euler}, which generalizes the Random Connection Model in a natural way and includes several models used in the literature as special cases.
We focus on the marked stationary case with vertices in $\R^d\times\AA$, where the mark space $\AA$ is an arbitrary Borel space.
We will derive a central limit theorem for an abstract class of functionals and show that many of the typical functionals considered in the study of simplicial complexes, such as Betti numbers, fall into this class.
As an important special case, we obtain a central limit theorem for Betti numbers in the Boolean model.
\end{abstract}

\begin{small}
\end{small}

\begin{small}
\mscs{60F05 \and 60D05 \and 60B99}
\end{small}

\section{Introduction}

Random graphs provide an effective mathematical framework for modeling and analyzing complex systems built on pairwise interactions.
Applications of random graphs span a wide range of disciplines, including the spread of diseases in epidemiology, communication and search algorithms in computer science, the structure of social and economic networks, and critical phenomena in statistical physics.
While graph models effectively capture pairwise interactions, their extension into higher dimensions through random simplicial complexes enables a deeper exploration of group dynamics and multi-way relationships.
In recent years, random simplicial complexes have gained significant attention in mathematics, ranging from combinatorial models, such as extensions of the classic Erd\H{o}s-Rényi graph (\cite{Costa,Kahle.Topology}), to applications in topological data analysis (see for example \cite{Biscio}), where the Vietoris-Rips and the \v Cech complex are used to study the structure of spatial data and extract topological features.

Betti numbers are among the most fundamental quantities studied in simplicial complexes, providing essential information about the topological structure of a space, such as its connectedness, loops, and voids. In topological data analysis, the Betti numbers play a crucial role in extracting topological features from data, offering a powerful way to quantify the underlying shape of datasets.
CLTs for the Betti numbers of the \v Cech and the Vietoris-Rips complex are, for example, provided in \cite{Kahle.LimitTheorems,Yogeshwaran}.
So-called soft random simplicial complexes, generalizations of the \v Cech and the Vietoris-Rips complex with additional randomness, have been introduced in \cite{Candela}.
Asymptotic results for Betti numbers of soft random simplicial complexes can be found in \cite{Candela.new}.

Formally, a \textit{simplicial complex} is a collection of non-empty finite subsets, known as simplices, of a vertex set $V$, with the condition that for each simplex, all of its faces (non-empty subsets) are also included.
A simplex containing $n+1$ vertices is called an $n$-simplex, and $n$ is also referred to as the dimension of the simplex.
The dimension of a simplicial complex is defined as the supremum of the dimensions of its simplices and can be infinite.
In this sense, graphs can be understood as 1-dimensional simplicial complexes.
The $j$-skeleton of a simplicial complex refers to the subcomplex formed by restricting the complex to simplices of dimension at most $j$.

In this paper we present a model extending the \textit{Random Connection Model (RCM)} to a random simplicial complex, as first proposed in \cite{Pabst.Euler}.
To this end, we fix a natural number $\alpha$, which determines the maximal dimension of the simplicial complex. 
In the special case $\alpha=1$, the model reduces to a random graph and coincides with the classical RCM, which has been studied as a random graph model for several decades (see, e.g., \cite{Caicedo,Can,Dickson,Last.RCM}) since it was introduced by Penrose (\cite{Penrose.RCM}) in 1991. 
In particular, all results obtained in this paper apply to random graphs as a special case of the present framework.
Although the model can be defined and studied in a more general setting (as in \cite{Pabst.Euler}), the methods employed in this paper rely on certain properties of Euclidean space.
Nevertheless, it turns out that we are not limited to $\R^d$ alone; we can, in fact, consider the product space $\R^d \times \AA$, where $(\AA, \mathcal{T})$ is an arbitrary Borel space equipped with a probability measure $\Theta$.

\begin{figure}[t!]
  \captionsetup{ labelfont = {bf}, format = plain }
  \centering
  \subfloat[][]{\includegraphics[width=0.48\linewidth]{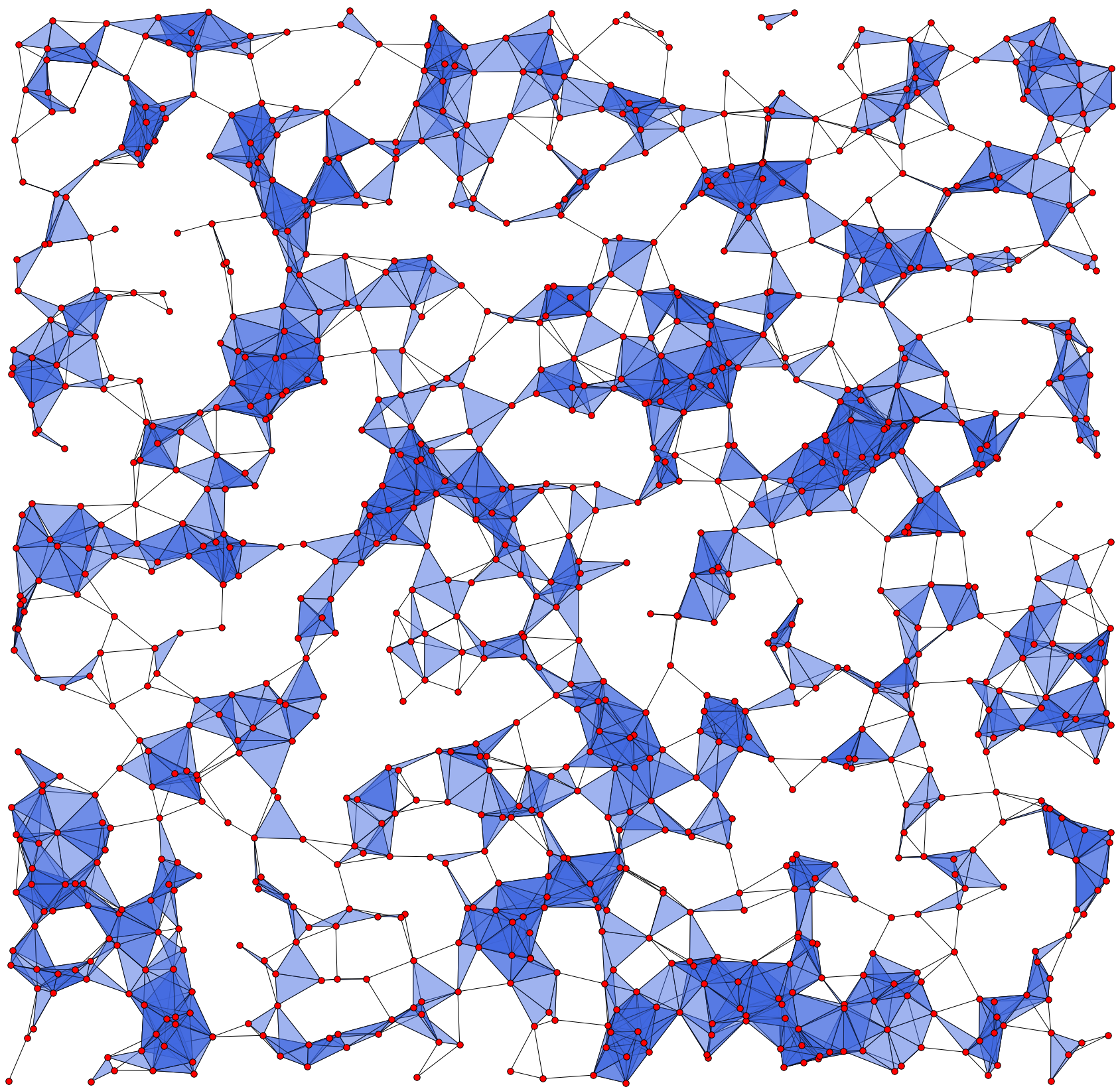}}
  \quad
  \subfloat[][]{\includegraphics[width=0.48\linewidth]{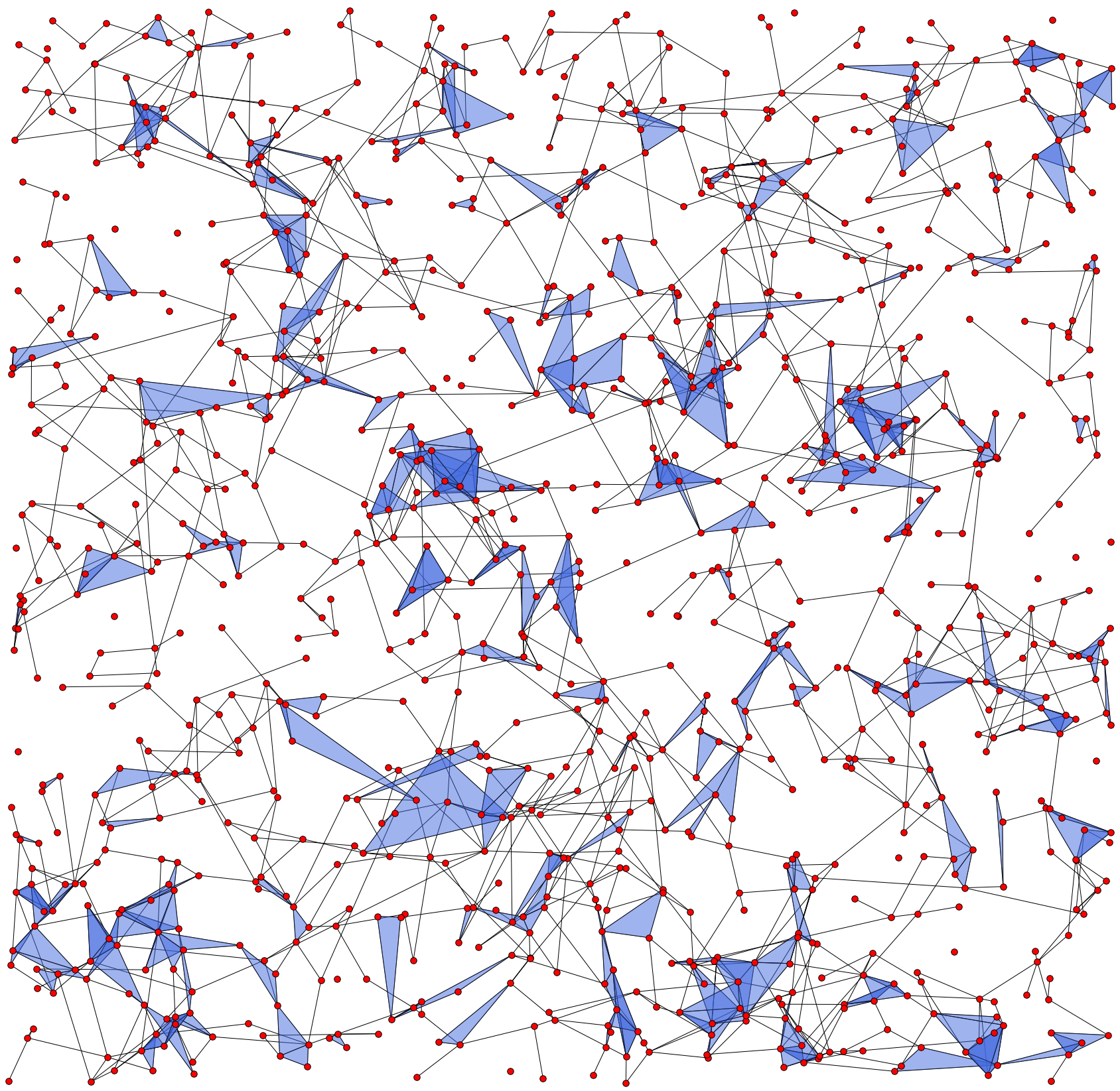}}
  \caption{(a) A realization of the model with the connection functions from (\ref{geometrisch}) and (b) a realization of the model with the connection functions from (\ref{exponentiell}) based on the same realization of a Poisson process}
  \label{fig:BspModell}
\end{figure}

The vertex set of the random simplicial complex will be the set of points of a Poisson process $\Phi$ on $\R^d\times\AA$ with intensity measure $\gamma\lambda_d\otimes\Theta$, where $\gamma>0$ is the intensity of the model and $\lambda_d$ denotes the Lebesgue measure on $\R^d$.
For all $j\in\{1,\dots,\alpha\}$ we fix a measurable, symmetric function $\varphi_j:(\R^d\times\AA)^{j+1}\rightarrow [0,1]$ which we assume to be additionally translation-invariant in the sense of (\ref{translation_invariant}).
We call the functions $\varphi_1,\dots,\varphi_{\alpha}$ the connection functions of the model and define the random simplicial complex $\Delta$ in the following way.
\begin{itemize}
\item[(0)] The vertex set of $\Delta$ is the set of points of $\Phi$.
\item[(1)] For each pair of points $x,y\in\Phi$ add the edge $\{x,y\}$ independently with probability $\varphi_1(x,y)$ to the complex $\Delta$.
\item[(2)] For each triple $x,y,z\in\Phi$, all of whose edges have been added to the complex, add the triangle $\{x,y,z\}$ independently with probability $\varphi_2(x,y,z)$ to the complex $\Delta$.
\item[$\vdots$]
\item[($\alpha$)] For $\alpha+1$ points $x_{i_0},\dots,x_{i_\alpha}\in\Phi$, all of whose subsimplices have been added to the complex, add the simplex $\{x_{i_0},\dots,x_{i_\alpha}\}$ independently with probability $\varphi_\alpha(x_{i_0},\dots,x_{i_\alpha})$ to the complex $\Delta$.     
\end{itemize}
We refer to the introduced framework as the \textit{marked stationary} case, which is consistent with the terminology used in the literature for the RCM.
We begin by illustrating the model through two examples.
As two-dimensional simplicial complexes are especially illustrative, we choose $\alpha = 2$.
In our first example model, we choose the connection functions
\begin{align}\label{geometrisch}
\varphi_1(x,y) \,=\, \1 \big\{ \Vert x-y \Vert \leq r \big\}, \qquad \varphi_2(x,y,z) \,\equiv\, p
\end{align}
for $r>0$ and $p\in [0,1]$.
The model is built in two steps: first, a geometric graph with parameter $r$ is formed; second, each potential triangle is added with fixed probability $p$.
For the second example, we consider the connection functions
\begin{align}\label{exponentiell}
\varphi_1(x,y) \,=\, \exp\big(-\theta_1 \norm{x-y}\big), \qquad \varphi_2(x,y,z) \,=\, \exp\big(-\theta_2 V(x,y,z)\big)
\end{align}
for $\theta_1,\theta_2>0$, where $V(x,y,z)$ denotes the (two-dimensional) volume of the triangle spanned by $x,y,z$.
Figure~\ref{fig:BspModell} shows realizations of both models in dimension $d=2$, within the unit square $[0,1]^2$.
Both are based on the same realization of a Poisson process, using $r=p=\frac{1}{2}$ for (a), and $\theta_1 = 4$, $\theta_2 = \frac{1}{100}$ for (b).

The connection functions can also be chosen such that $\Delta$ coincides with the $\alpha$-skeleton of the \v Cech complex, the Vietoris–Rips complex, or the complexes considered in \cite{Candela}.
However, for functionals for which there exists $j\in\N$ such that they depend only on the $j$-skeleton of a simplicial complex, this does not impose any restriction, since one can simply choose $\alpha\geq j$.
This applies to all concrete example functionals considered in this paper, except for the Euler characteristic.
For example, the $p$-th Betti number of a simplicial complex is already determined by its $(p+1)$-skeleton.
In Section~\ref{Sec:Euler}, we will see that the CLT can nevertheless be applied to the Euler characteristic of both the \v Cech and the Vietoris–Rips complex.

Our main result is a CLT that applies to various functionals, including the \textit{Betti numbers}.
The proof builds on the methods developed in \cite{Can}, where the RCM is studied as a random graph in the unmarked stationary setting.
Note that the unmarked stationary case discussed in \cite{Can} arises as a special case of the marked framework presented here by taking the mark space $\AA$ to consist of a single element.
Consequently, the CLT we prove is novel in two aspects: it applies to simplicial complexes with dimensions greater than one and also addresses the marked case, even specifically for graphs.
In \cite{Can}, the CLT is applied to the Betti numbers of the clique complex of the RCM.
The clique complex of a graph is the complex formed by adding all simplices of dimension greater than 1 whose edges are already present in the graph.
To the best of our knowledge, the result in \cite{Can} represents the pioneering CLT addressing higher Betti numbers (excluding the zeroth Betti number) within the context of the RCM.
This result can be recovered as a special case of the present framework by choosing all connection functions except for $\varphi_1$ to be constant and equal to 1 in the unmarked case.
In contrast, our results apply directly to the first Betti number $\beta_1$ of the underlying random graph (corresponding to the case $\alpha=1$).
In particular, this yields a CLT for $\beta_1$ of the random graph, which is new even in the unmarked case.
A CLT for the number of components, i.e., the zeroth Betti number, can be found in \cite{Last.RCM} (Theorem~9.1).
The paper \cite{Penrose.IRG} studies the asymptotic behavior of various functionals in the RCM on a general space.
These include the number of vertices of a fixed degree, as well as the number of induced subgraphs or connected components of a given order.

An important application of our result refers to the \textit{Boolean model}, one of the flagship models of stochastic geometry.
In \cite{Hug}, a CLT (Theorem~9.1) for a large class of functionals in the Boolean model is presented.
However, this does not include the Betti numbers, as the functionals considered there are additive.
Indeed, as far as we know, CLTs for Betti numbers in the Boolean model have been known only in a very restrictive special case, specifically when balls with identical radii are used as grains.
This special case corresponds to the \v Cech complex, and a corresponding result can be found, for example, in \cite{Yogeshwaran} (Theorem~4.7).
For more details on the Boolean model, we refer to Section~\ref{Sec:BooleanModel}.
In this paper, results and methods derived from the PhD thesis \cite{Pabst.Thesis} are presented.

The main results of the paper are stated in Section~\ref{Sec:MainResults}.
After introducing the required preliminaries in Section~\ref{Sec:Pre} and providing a precise definition of the model in Section~\ref{Sec:Model_Introduction}, Section~\ref{Sec:CLT} establishes and proves the CLT for weakly stabilizing functionals.
Section~\ref{Sec:weaklyStabFct} studies the class of weakly stabilizing functionals.
We show that many typical functionals arising in the context of graphs and simplicial complexes are weakly stabilizing under suitable integrability assumptions.
Section~\ref{Sec:Variance} addresses the question of when the asymptotic variance in the CLT is strictly positive.
Finally, Section~\ref{Sec:BooleanModel} illustrates how to derive a CLT for Betti numbers within the Boolean model.

\section{Main results} \label{Sec:MainResults}

Our main result is a CLT for so-called weakly stabilizing functionals according to the upcoming definition.
We denote by $\Delta$ the random simplicial complex described in the introduction.
A precise definition of the model is given in Section~\ref{Sec:Model_Introduction}.
For a measurable set $W\subset\R^d$, let $\Delta_W$ denote the restriction of $\Delta$ to simplices with vertices in $W\times\AA$.
A key tool for deriving a CLT and defining weakly stabilizing functionals is a suitable difference operator, which measures the influence of adding an additional vertex to the model.

To define the difference operator, we work with an augmented Poisson process $\Psi$ on $\R^d\times\AA\times[0,1]\times\MM$, in which each point carries additional auxiliary marks taking values in $[0,1]\times\MM$.
Here $\MM$ denotes a countable product of copies of the unit interval.
The auxiliary marks in $[0,1]\times\MM$ are used to determine which simplices are included in the complex, but they are not part of the vertex labels of the resulting simplicial complex.
The full construction can be found in Section~\ref{Sec:Model_Definition}.
For $(x,a,t,u)\in\R^d\times\AA\times[0,1]\times\MM$, let $\Delta^{(x,a,t,u)}$ denote the simplicial complex obtained from the augmented point process $\Psi+\delta_{(x,a,t,u)}$.
In particular, the vertices of $\Delta^{(x,a,t,u)}$ lie in $\R^d\times\AA$.
For a measurable, bounded set $W\subseteq\R^d$ and a real-valued function $f$ defined for finite simplicial complexes, we define the difference operator
\begin{align}
\Lambda_{(x,a,t,u)}f(\Delta_{W}) \,:=\, f\big(\Delta_{W}^{(x,a,t,u)}\big) - f\big(\Delta_{W\setminus\{x\}}^{(x,a,t,u)}\big).
\end{align}
By construction, the complex $\Delta_{W\setminus\{x\}}^{(x,a,t,u)}$ is obtained from $\Delta_{W}^{(x,a,t,u)}$ by removing the additional vertex $(x,a)$ and is therefore a subcomplex.
The operator $\Lambda_{(x,a,t,u)}f$ thus quantifies the effect of adding the
vertex $(x,a)$ on the value of the functional $f$.
In the main results below, this operator is evaluated at a typical added point of the form $(\0,V,1,U)$, where $V\sim\Theta$ and $U\sim\Q$ are independent random elements, in particular independent of $\Psi$.
Here $\Theta$ denotes the mark distribution on $\AA$, and $\Q$ is a countable
product of the uniform distribution on $[0,1]$.
Due to the translation-invariant situation, the Euclidean component of the added point does not matter, and we choose the origin $\0$ as usual.

Furthermore, for a simplicial complex $K$ with vertices in $\R^d\times\AA$, we denote by $K+t$, $t\in\R^d$, the simplicial complex obtained from $K$ by translating all Euclidean components of vertices of $K$ by $t$.
That is, each simplex $\{(x_0,a_0),\dots,(x_j,a_j)\}$, $j\in\N_0$, is replaced by the translated simplex $\{(x_0+t,a_0),\dots,(x_j+t,a_j)\}$.
Moreover, let $\W$ be the set of all (half-open) cubes in $\R^d$.
Finally, when we refer to a sequence of cubes $(W_n)_{n\in\N}$ tending to $\R^d$, briefly denoted by $W_n\rightarrow\R^d$, we mean that for every bounded set $B\subset\R^d$, there exists an $N\in\N$ such that $B\subseteq W_n$ holds for all $n\geq N$.

\begin{Def} [Weakly stabilizing functionals] \label{Def:weaklyStabFunctionals}
Let $f$ be a measurable real-valued function defined for finite simplicial complexes with vertices in $\R^d\times\AA$.
\begin{enumerate}
\item[$\mathrm{(i)}$] The function $f$ is called translation invariant if $f(K)=f(K+t)$ holds for all $t\in\R^d$ and all finite simplicial complexes $K$ with vertices in $\R^d\times\AA$.
\item[$\mathrm{(ii)}$] The function $f$ is called weakly stabilizing if it is translation invariant and there exists a real-valued random variable $Z$ such that for all sequences of cubes $(W_n)_{n\in\N}$ with $W_n\rightarrow\R^d$, the convergence
\begin{align} \label{Convergence_weaklyStab}
\Lambda_{(\0,V,1,U)}f\big(\Delta_{W_n}\big) \;\xrightarrow{\P}\; Z
\end{align}
holds as $n\rightarrow\infty$.
We call $Z$ the random variable associated with $f$.
\item[$\mathrm{(iii)}$] We say $f$ satisfies the moment condition if there is some $\varepsilon>0$ with
\begin{align} \label{MomentCondition}
\sup\limits_{\substack{W\in\W, \\ \0\in W}} \E{ \big|\Lambda_{(\0,V,1,U)} f(\Delta_W)\big|^{2+\varepsilon} } \,<\, \infty.
\end{align}
\end{enumerate}
\end{Def}

Observe that for a weakly stabilizing functional $f$, the random variable $Z$ associated to $f$ is almost everywhere uniquely determined, which justifies the usage of the definite article.
If $f$ also satisfies the moment condition (\ref{MomentCondition}), the convergence in (\ref{Convergence_weaklyStab}) even holds in $L^2$.
Before we address concrete example functionals in Section~\ref{Sec:weaklyStabFct}, let us first note that the class of weakly stabilizing functionals satisfying the moment condition is non-empty, as all constant functionals belong to this class.
Moreover, the class is closed under taking linear combinations.
The moment condition for such a linear combination follows by applying the Minkowski inequality.

The following lemma establishes the existence of an asymptotic variance and forms a crucial step in proving a CLT for weakly stabilizing functionals.
It also provides a criterion for the asymptotic variance to be strictly positive, which serves as the starting point for a later analysis.

\begin{Lem} \label{Lem:asymptoticVariance}
Let $f$ be a weakly stabilizing functional that satisfies the moment condition (\ref{MomentCondition}).
There exists $\sigma^2\in[0,\infty)$ such that for every sequence of cubes $(W_n)_{n\in\N}$ with $W_n\rightarrow\R^d$, the convergence
\begin{align*}
\frac{\V{\fWn}}{|W_n|} \;\rightarrow\; \sigma^2 \qquad \text{as } n\rightarrow\infty
\end{align*}
holds.
If the random variable $Z$ associated with $f$ from Definition~\ref{Def:weaklyStabFunctionals} satisfies $\P( Z\neq 0 )>0$, then $\sigma^2>0$.
\end{Lem}

An explicit representation of $\sigma^2$ is given in Lemma~\ref{Lem:VarRepresentation}.
The main result is the following CLT for weakly stabilizing functionals.
A multivariate version is given in Section~\ref{Sec:CLT}.

\begin{Th} [CLT for weakly stabilizing functionals] \label{Th:ZGWS3}
Let $f$ be a weakly stabilizing functional satisfying the moment condition (\ref{MomentCondition}).  
Then for any sequence of cubes $(W_n)_{n\in\N}$ with $W_n\rightarrow\R^d$, we have  
\begin{align*}
\frac{f\big(\Delta_{W_n}\big) - \E{ f\big(\Delta_{W_n}\big) }}{\sqrt{|W_n|}} \;\xrightarrow{d}\; N_{\sigma^2} \qquad \text{as } n \to \infty,
\end{align*}
where $\sigma^2$ is the limit from Lemma~\ref{Lem:asymptoticVariance}.  
If $\P( Z \neq 0 ) > 0$ for the random variable $Z$ associated with $f$, then the asymptotic variance $\sigma^2$ is positive.
\end{Th}

A natural question is whether one can formulate concrete conditions on the
function $f$ and the connection functions that ensure positivity of the
asymptotic variance.
To this end, we will derive a corollary that provides a condition under which the asymptotic variance is indeed strictly positive for weakly additive functionals.
Here we call a functional $f$ weakly additive if, for a finite simplicial complex $K$ with connected components $Z_1,\dots,Z_m$, the equality
\begin{align*}
    f(K) \,=\, \sum_{i=1}^m f(Z_i)
\end{align*}
always holds.
This property was formulated for the Betti numbers in Lemma~\ref{Lem:Betti_properties} (iii).
For a finite set of points $(x_1,a_1),\dots,(x_k,a_k)\in\R^d\times\AA$, we denote by $\Delta\big((x_1,a_1),\dots,(x_k,a_k)\big)$ the random simplicial complex we obtain by applying the construction of $\Delta$ (in terms of the multi-staged random experiment described in the introduction) to the vertex set $\{(x_1,a_1),\dots,(x_k,a_k)\}$.

\begin{Kor} \label{Kor:Positivity}
Let $f$ be a function defined on finite simplicial complexes, which is weakly additive, weakly stabilizing, invariant under isomorphism and satisfies the moment condition (\ref{MomentCondition}).
Assume that there exists a finite connected simplicial complex $K$ on $r+1$ vertices, $r \in \mathbb{N}_0$, with $\dim(K) \leq \alpha$, such that $|f(K) - f(L)| > 0$ for every induced subcomplex $L$ of $K$ on $r$ vertices.
Additionally, suppose that one of the following two conditions is fulfilled.
\begin{enumerate}
    \item[$\mathrm{(i)}$] The connection functions are such that $\varphi_j\in (0,1)$ almost everywhere for each $j\leq\min\{\alpha,\dim(K)+1\}$.
    \item[$\mathrm{(ii)}$] The simplicial complex $K$ satisfies
    \begin{align*}
        \int_{\left( \R^d \right)^{r}} \int_{\AA^{r+1}} \P\left( \Delta\big((\0,a_0),(x_1,a_1),\dots,(x_r,a_r)\big) \cong K \right) \; \Theta^{r+1} (\d a) \; \d x \,>\, 0,
    \end{align*}
    where we abbreviate $\d a = \d {(a_0,\dots,a_r)}$ and $\d x = \d {(x_1,\dots,x_r)}$.
\end{enumerate}
Then the asymptotic variance $\sigma^2$ from Theorem~\ref{Th:ZGWS3} applied to $f$ is strictly positive.
\end{Kor}

\section{Preliminaries} \label{Sec:Pre}

This section provides some mathematical background used in this paper, beginning with the notion of a Borel space.
A Borel space $(\AA,\mathcal{T})$ is a measurable space such that there is a bijective measurable map $\phi:\AA\rightarrow U$ into some measurable subset $U$ of the unit interval $[0,1]$ with measurable inverse.
We assume that all random elements appearing in this work live on some probability space $(\Omega,\A,\P)$.

\subsection{Poisson processes}

To introduce Poisson processes on $\R^d\times\AA$ let $\ND$ be the set of all locally finite counting measures on $\R^d\times\AA$.
Here a counting measure $\eta$ on $\R^d\times\AA$ is called locally finite, when $\eta(B\times\AA)<\infty$ holds for all bounded Borel sets $B$.
We equip $\ND$ with the smallest $\sigma$-field such that the mappings $\ND\rightarrow [0,\infty]$, $\eta\mapsto\eta(B)$ are measurable for all measurable sets $B\subseteq\R^d\times\AA$.
A random element $\Phi$ in $\ND$ is called a Poisson process on $\R^d\times\AA$ with intensity measure $\gamma\lambda_d\otimes\Theta$, if the following two properties are satisfied.
\begin{enumerate}
    \item[(a)] For every measurable set $B\subseteq\R^d\times\AA$ the random variable $\Phi(B)$ follows a Poisson distribution with parameter $(\gamma\lambda_d\otimes\Theta)(B)$.
    \item[(b)] For pairwise disjoint measurable sets $B_1,\dots,B_m\subseteq\R^d\times\AA$, $m\in\N$, the random variables $\Phi(B_1),\dots,\Phi(B_m)$ are stochastically independent.
\end{enumerate}
Theorem~3.6 in \cite{Last.Lectures} ensures that a Poisson process in this setting always exists and its distribution is uniquely determined by $\gamma$ and $\Theta$ (Proposition~3.2 in \cite{Last.Lectures}).
Throughout this paper $\Phi$ will always be a Poisson process on $\R^d\times\AA$ with intensity measure $\gamma\lambda_d\otimes\Theta$.
For a measurable set $B\subseteq\R^d$ we will denote by $\Phi_B$ the restriction of $\Phi$ to $B\times\AA$.

\subsection{Simplicial complexes and Betti numbers} \label{Sec:Definition_Betti}

A simplicial complex is a family of nonempty, finite subsets of a vertex set, which is closed under taking subsets.
An element $\sigma$ of a simplicial complex is called a $j$-simplex if it has $j+1$ elements and $j$ is called the dimension of the simplex.
We will denote a $j$-simplex with vertices $v_0,\dots,v_j$ by $[v_0,\dots,v_j]$.
Two simplicial complexes $K,L$ with vertex sets $V_K,V_L$ are called isomorphic (denoted by $K\cong L$) if there is a bijective map $f:V_K\rightarrow V_L$ having the property
\begin{align*}
\sigma\in K \quad \Longleftrightarrow \quad f(\sigma)\in L.
\end{align*}
For a vertex $v$, that is a 0-simplex, of a simplicial complex $K$ and $n\in\N$, we call
\begin{align} \label{SimplexDegree}
    \deg_n(v,K) \,:=\, \|\{ \sigma\in K \,\|\, v\in\sigma, \dim(\sigma)=n \}\|
\end{align}
the $n$-simplex degree of the vertex $v$ in $K$.
The $j$-skeleton of a simplicial complex is its restriction
\begin{align*}
S_j(K) \,:=\, \{ \sigma\in K \,\|\, \dim(\sigma)\leq j \}
\end{align*}
to simplices of dimension not bigger than $j$ and is again a simplicial complex.
Furthermore we denote by $f_i(K)$ the number of $i$-simplices of a simplicial complex $K$.
Throughout this paper we will only consider finite simplicial complexes that are complexes containing only a finite number of simplices.
While we will derive a CLT for a broad range of functionals of simplicial complexes, the Betti numbers are of particular interest to us.
The Betti numbers are intricate quantities, and their definition involves several preparatory steps.
For this purpose, we fix a finite simplicial complex $K$  and denote the field containing only the two elements 0 and 1 by $\Z_2$.
For $p\in\N_0$, let $\sigma_1,\dots,\sigma_k$ be the $p$-simplices of $K$.
We define the $p$-th chain group $C_p(K)$ of $K$ as the vector space over $\Z_2$, for which $\{\sigma_1,\dots,\sigma_k\}$ is a basis, that is
\begin{align*}
    C_p(K) \,:=\, \Big\{ \sum_{i=1}^k a_i\sigma_i \,\big|\, a_i\in \Z_2 \Big\}.
\end{align*}
Due to the choice of $\Z_2$ as the space of coefficients (which is common in applications), $C_p(K)$ can be identified with the power set of $\{\sigma_1,\dots,\sigma_k\}$.
For $p\in\N$ define a linear map  $\partial_p:C_p(K)\rightarrow C_{p-1}(K)$ by linear extension of
\begin{align*}
\partial_p [v_0,\dots,v_p] \,:=\, \sum_{j=0}^p [v_0,\dots,v_{j-1},v_{j+1},\dots,v_p].
\end{align*}
With the additional definitions $C_{-1}(K):=\{0\}$ and $\partial_0\equiv 0$ we obtain a so-called chain complex

\begin{center}
\begin{tikzpicture}[xscale=0.8,yscale=0.8]
\node (A) at (0,0) {$\dots$};
\node (B) at (2.5,0) {$C_{p+1}(K)$};
\node (C) at (5.5,0) {$C_{p}(K)$};
\node (D) at (8.5,0) {$C_{p-1}(K)$};
\node (E) at (11.2,0) {$\dots$};
\node (F) at (13.7,0) {$C_{0}(K)$};
\node (G) at (16.3,0) {$\{0\},$};
\draw[->, thin] (A) to (B);
\draw[->, thin] (B) to node[above] {\footnotesize $\partial_{p+1}$} (C);
\draw[->, thin] (C) to node[above] {\footnotesize $\partial_{p}$} (D);
\draw[->, thin] (D) to node[above] {\footnotesize $\partial_{p-1}$} (E);
\draw[->, thin] (E) to node[above] {\footnotesize $\partial_{1}$} (F);
\draw[->, thin] (F) to node[above] {\footnotesize $\partial_{0}$} (G);
\end{tikzpicture}
\end{center}
satisfying $\partial_p\circ\partial_{p+1}\equiv 0$.
The last condition is equivalent to
\begin{align*}
B_p(K) \,:=\, \mathrm{image}(\partial_{p+1}) \,\subseteq\, \mathrm{kernel}(\partial_p) \,=:\, Z_p(K).
\end{align*}
The $p$-th Betti number $\beta_p(K)$ of $K$ is defined by
\begin{align*}
\beta_p(K) \,:=\, \dim(Z_p(K)) - \dim(B_p(K)),
\end{align*}
which is the dimension of the quotient vector space $Z_p(K)/ B_p(K)$.
We summarize the most important properties of the Betti numbers in the following lemma.

\begin{Lem} \label{Lem:Betti_properties}
    Let $p\in\N_0$ and $K$ be a finite simplicial complex whose connected components are given by $Z_1,\dots,Z_m$.
    Furthermore, let $L\subseteq K$ be a subcomplex of $K$.
    Then the following statements hold.
    \begin{enumerate}
        \item[$\mathrm{(i)}$] $\beta_p(K)=0$\; if $p>\dim(K)$.
        \item[$\mathrm{(ii)}$] $\beta_p(K)=\beta_p(S_{p+1}(K))$.
        \item[$\mathrm{(iii)}$] $\beta_p(K)=\sum_{i=1}^m\beta_p(Z_i)$.
        \item[$\mathrm{(iv)}$] $|\beta_p(K)-\beta_p(L)|\,\leq\, (f_p(K)-f_p(L))+(f_{p+1}(K)-f_{p+1}(L))$.
    \end{enumerate}
\end{Lem}

The first two properties follow immediately from the definition of the Betti numbers.
A proof of the third one, which sometimes is called weakly additive, can be found in \cite{Hatcher} (Proposition~2.6).
Property (iv) of Lemma~\ref{Lem:Betti_properties} is crucial to obtain a CLT for the Betti numbers, and is proved, for example, in \cite{Yogeshwaran} (Lemma 2.2).

\section{The RCM for higher-dimensional simplicial complexes} \label{Sec:Model_Introduction}

In this section, we provide a precise definition of the random simplicial complex $\Delta$. 
We also establish variance results related to the difference operator introduced in Section~\ref{Sec:MainResults}.

\subsection{Definition of the model} \label{Sec:Model_Definition}

First of all, we fix some natural number $\alpha\in\N$, which will be the maximum dimension of the simplicial complex $\Delta$.
Furthermore, for each $j\in\{1,\dots,\alpha\}$, we choose a measurable, symmetric and translation-invariant function $\varphi_j:(\R^d\times\AA)^{j+1}\rightarrow [0,1]$, where the latter means that
\begin{align} \label{translation_invariant}
    \varphi_j\big((x_0+t,a_0),\dots,(x_j+t,a_j)\big) \,&=\, \varphi_j\big((x_0,a_0),\dots,(x_j,a_j)\big)
\end{align}
holds for all $x_0,\dots,x_j,t\in\R^d$ and $a_0,\dots,a_j\in\AA$.
We call the functions $\varphi_1,\dots,\varphi_\alpha$ the connection functions of the model.
For the definition of $\Delta$ we need an even bigger space than $\R^d\times\AA$, where the additional randomness to decide which simplices are part of the complex comes from.
We write
\begin{align*}
M^{(j)} \,:=\, [0,1]^{\N^j} \,=\, \big\{ (u_{z})_{z\in\N^j} \,\|\, u_{z}\in [0,1] \text{ for all } z\in\N^j \big\}
\end{align*}
for the space of $j$-times indexed sequences with members in $[0,1]$, equipped with the product $\sigma$-field of the Borel $\sigma$-fields on $[0,1]$.
We consider the space $\MM := M^{(2)} \times \dots \times M^{(2\alpha)}$ together with the product $\sigma$-field and the probability distribution $\Q:=\otimes_{j=1}^\alpha\otimes_{z\in\N^{2j}}\mathcal{U}([0,1])$, where $\mathcal{U}([0,1])$ is the uniform distribution on $[0,1]$.
A random element $U$ in $\MM$ with distribution $\Q$ is therefore a countable sequence of i.i.d. uniformly distributed random variables.
For $u\in\MM$, $j\in\{1,\dots,\alpha\}$ and $z=(k_1,l_1,\dots,k_j,l_j)\in\N^{2j}$ we write $u_{k_1,l_1,\dots,k_j,l_j}$ for the component of $u$ associated with $z$.

Now let $\Psi$ be a Poisson process on $\R^d\times\AA\times [0,1]\times\MM$ with intensity measure $\gamma\lambda_d\otimes\Theta\otimes\mathcal{U}([0,1])\otimes\Q$ and $\Phi$ its projection on $\R^d\times\AA$, which is a Poisson process with intensity measure $\gamma\lambda_d\otimes\Theta$.
Observe that, since a countable product of Borel spaces is again a Borel space, $\Psi$ is again a Poisson process on a product of $\R^d$ with some Borel space.
For a point $(x,a,t,u)\in\Psi$, the coordinate $t\in[0,1]$ is used to impose an ordering on the points of the process, while the components of $u\in\MM$ will be used to decide which simplices are contained in $\Delta$.

We define the vertex set of $\Delta$ as the set of points of $\Phi$.
To determine which higher-dimensional simplices belong to the complex $\Delta$, we associate to each potential simplex $\sigma$ a coordinate $u_\sigma$ taken from one of the marks $u\in\MM$ of the points of $\Psi$.
For this purpose, we introduce a suitable ordering of the points of $\Psi$.
We assign to each point in $\Psi$ a unique pair coordinates $(k,l)\in\N_0^2$ as follows.
Consider the partition $(Q_i)_{i\in\N_0}$ of $\R^d$ into cubes $Q_i := [0,1)^d + z_i$, where $z_0,z_1,\dots$ is a fixed enumeration of $\Z^d$ with $z_0=\0$ (origin).
For a point $(x,a,t,u)\in\Psi$, the index $k$ is determined uniquely by the condition $x\in Q_k$, and $l$ is such that the point is the $l$-th smallest point of $\Psi$ in $Q_k$ with respect to the ordering defined by
\begin{align} \label{ordering_points}
(x,a,t,u) \,<\, (y,b,s,v) \qquad \text{if } t < s.
\end{align}
With probability one, each point of $\Psi$ is assigned a unique pair of coordinates $(k,l)\in\N_0^2$, which identifies the point within the process.

Fix points $(x_0,a_0,t_0,u^{(0)}),\dots,(x_j,a_j,t_j,u^{(j)})\in \Psi$, $j\in\{1,\dots,\alpha\}$, with $t_0<\dots<t_j$ and consider the potential simplex $\sigma:=[(x_0,a_0),\dots,(x_j,a_j)]$ of $\Delta$.
To this simplex, we assign the coordinate
\begin{align*}
u(\sigma) \,:=\, u^{(j)}_{k_0,l_0,\dots,k_{j-1},l_{j-1}},
\end{align*}
where $(k_i,l_i)$ denotes the coordinates of the point $(x_i,a_i,t_i,u^{(i)})$ within $\Psi$.
We define the rule
\begin{align*}
    \sigma\in\Delta \enskip :\Longleftrightarrow \enskip u(\rho)\leq \varphi_{|\rho|-1}(x_\rho) \enskip \text{for all } \rho\subseteq\sigma \text{ with } |\rho|\geq 2,
\end{align*}
where $\varphi_{|\rho|-1}(x_\rho)$ is the value of $\varphi_{|\rho|-1}$, when inserting the elements of $\rho$.
This decision rule uniquely determines the simplicial complex $\Delta$.
Observe that translating the Euclidean components of all points of $\Psi$ by the same vector while keeping their marks fixed results in another simplicial complex (with the same distribution), which, on the level of realizations, need not be isomorphic to the original one.
For instance, two points that were connected by an edge before the shift may no longer be connected afterwards.

For $W\subseteq\R^d$, we write $\Delta_W$ for the restriction of $\Delta$
to simplices with vertices in $W\times\AA$, and $\Psi_W$ for the restriction
of $\Psi$ to $W\times\AA\times[0,1]\times\MM$.
Since the connection functions and the distribution of the Poisson process $\Psi$ are translation invariant with respect to the Euclidean component, it follows that for every $t\in\mathbb{R}^d$, the random variables $f(\Delta_W)$ and $f(\Delta_{W+t})$ have the same distribution whenever $f$ is translation invariant in the sense of Definition~\ref{Def:weaklyStabFunctionals}.

This construction has another important property, which we briefly discuss and which also motivates its particular form.
Whenever $W$ is a union of cubes from the decomposition $(Q_i)_{i\in\N_0}$, the subcomplex $\Delta_W$ is fully determined by $\Psi_W$, as this holds for the coordinates of the points of $\Psi_W$ within $\Psi$.
Consequently, for two disjoint sets $W_1,W_2$, which are unions of cubes from the decomposition $(Q_i)_{i\in\N_0}$, the subcomplexes $\Delta_{W_1},\Delta_{W_2}$ are independent since the processes $\Psi_{W_1},\Psi_{W_2}$ are independent.
This property will be crucial for the derivation of a CLT.

The idea of introducing an additional $[0,1]$-component in the mark space to order the points of the process already appears in \cite{Last.RCM,Penrose.IRG}, among other constructions.
We refer in particular to the approach developed in Section~5 of \cite{Last.RCM}.
That approach also relies on a decomposition of the underlying space and an implicit introduction of coordinates for the points, analogous to those used here.
However, in \cite{Last.RCM} additional assumptions on the decomposition are imposed which are not required in the present setting.  
This framework was subsequently adapted in \cite{Can} to the case of $\mathbb{R}^d$ using a decomposition into cubes.  
In all of these works, the objects under consideration are random graphs, which lead to a simpler overall construction and notation compared to the random simplicial complex studied here.
In the special case $\alpha=1$, the present construction reduces to a random graph.

\subsection{Difference operator and variance results}

In what follows, we denote by $\Delta^{(x,a,t,u)}$, $(x,a,t,u)\in\R^d\times\AA\times [0,1]\times\MM$, the simplicial complex constructed from the point process $\Psi+\delta_{(x,a,t,u)}$ using the previously introduced procedure from Section~\ref{Sec:Model_Definition}.
Observe that for $t=1$, the complex $\Delta^{(x,a,1,u)}$ contains $\Delta$ as a subcomplex almost surely, since the added point does not change the coordinates of the other points.
Since this no longer holds for general $t\in [0,1]$, we will work with the difference operator introduced in Section~\ref{Sec:MainResults} and recall its definition.
For a measurable, bounded set $W\subseteq\R^d$ and a real-valued function $f$ defined for finite simplicial complexes, let
\begin{align} \label{Def:DiffOp}
\Lambda_{(x,a,t,u)}f(\Delta_{W}) \,:=\, f\big(\Delta_{W}^{(x,a,t,u)}\big) - f\big(\Delta_{W\setminus\{x\}}^{(x,a,t,u)}\big).
\end{align}
The complex $\Delta_{W\setminus\{x\}}^{(x,a,t,u)}$ is a subcomplex of $\Delta_{W}^{(x,a,t,u)}$ since it is obtained from the latter by removing the vertex $(x,a)$ together with all simplices containing it.
This is the crucial property of this operator, ensuring that it depends only on the resulting complex and the location of the added point and not on how the ordering is affected.
Therefore, it is more suitable for our analysis than the classical difference operator for Poisson functionals
\begin{align*}
    D_{(x,a,t,u)}f(\Delta_{W}) \,:=\, f\big(\Delta_{W}^{(x,a,t,u)}\big) - f(\Delta_{W})
\end{align*}
from the literature, where one compares the functional before and after adding a point.
Here the second complex is not a subcomplex of the first one, because the added point changes the ordering and thus the assignment of the i.i.d. uniformly distributed random variables to the simplices.

The idea of defining such a well-behaved operator as in (\ref{Def:DiffOp}) was already used in \cite{Last.RCM} and later in \cite{Can}.
Observe that in principle one could define something like the "space of simplicial complexes" on which $f$ lives, but for this paper this is not necessary and one can think of quantities like $f(\Delta_W)$ or $\Lambda_{(x,a,t,u)}f(\Delta_{W})$ as functions of $\Psi$.

Now consider independent typical marks $V\sim\Theta$ and $U\sim\mathbb{Q}$, also independent of $\Psi$.
Due to the translation-invariant framework introduced in Section~\ref{Sec:Model_Definition}, it follows that for translation-invariant functionals $f$ (cf. Definition~\ref{Def:weaklyStabFunctionals}),
\begin{align*}
	\Lambda_{(x,V,t,U)} f(\Delta_W) \,\stackrel{d}{=}\, \Lambda_{(\0,V,t,U)} f(\Delta_{W-x}).
\end{align*}
Furthermore, the mark $t$ only affects the ordering of the points.
Hence, choosing a different value for $t$ may change the simplicial complex on the level of realizations but not its distribution, i.e.
\begin{align*}
	\Lambda_{(\0,V,t,U)} f(\Delta_W) \,\stackrel{d}{=}\, \Lambda_{(\0,V,t^\prime,U)} f(\Delta_{W}) \qquad\text{for all } t,t^\prime\in [0,1].
\end{align*}

A crucial step in deriving CLTs for Poisson functionals is controlling their variance.
To this end, the theory of Poisson processes provides variance formulas that utilize classical difference operators.
However, it has been shown that these can be transferred to operators like those from (\ref{Def:DiffOp}).
In the following, we present two such results that will be applied in later sections.
The first one is a generalization of Theorem~5.1 in \cite{Last.RCM}.
To state it, we denote by $\Psi_t$ the restriction of $\Psi$ to $\R^d\times\AA\times [0,t)\times\MM$.

\begin{Lem} \label{Lem:exactVar}
    Let $W\subseteq\R^d$ be measurable and bounded and $f$ a function defined for finite simplicial complexes with vertices in $\R^d\times\AA$, such that $\E{f(\Delta_W)^2}<\infty$.
    Then
    \begin{align*}
\V{f(\Delta_W)} = \gamma\int_0^1 \int_{W} \int_\AA \int_\MM \E{ \E{ \Lambda_{(x,a,t,u)}f(\Delta_W) \;\big|\; \Psi_t}^2 } \Q \, (\d {u}) \; \Theta \, (\d {a})\; \d {x} \; \d t.
\end{align*}
\end{Lem}

\begin{proof}
Since $W$ is bounded, it intersects only finitely many cubes of the decomposition $(Q_i)_{i\in\N_0}$ used for the construction of $\Delta$.
Let $Q_W$ denote the union of those cubes intersecting $W$.
From the construction of the simplicial complex, it follows that the difference operator $D_{(x,a,t,u)}f(\Delta_W)$ vanishes whenever $x\notin Q_W$.
Therefore, Theorem~B.1 in \cite{Last.RCM} shows
\begin{align*}
\V{f(\Delta_W)} = \gamma\int_0^1 \int_{Q_W} \int_\AA \int_\MM \E{ \E{ D_{(x,a,t,u)}f(\Delta_W) \;\big|\; \Psi_t}^2 } \Q \, (\d {u}) \; \Theta \, (\d {a})\; \d {x} \; \d t.
\end{align*}
By construction, $\Delta_W\stackrel{d}{=}\Delta_{W\setminus\{x\}}^{(x,a,t,u)}$, and hence $\E{f\big(\Delta_{W\setminus\{x\}}^{(x,a,t,u)}\big)^2}<\infty$.
Observe that inserting $(x,a,t,u)$ can only affect simplices containing at least one point $(y,b,s,v)$ with $s>t$.
Since the presence of each simplex is determined by the corresponding mark in $\MM$ of its largest point with respect to the ordering \eqref{ordering_points}, the complexes $\Delta_W$ and $\Delta_{W\setminus\{x\}}^{(x,a,t,u)}$ have the same conditional distribution given $\Psi_t$.
Consequently,
\begin{align} \label{Eq:bedEW}
\E{ f(\Delta_W) \;\big|\; \Psi_t} \,=\, \E{ f\big(\Delta_{W\setminus\{x\}}^{(x,a,t,u)}\big) \;\big|\; \Psi_t} \qquad \P\text{-almost surely}.
\end{align}
Applying the Mecke formula to the Poisson process $\Psi$ and the function
\begin{align*}
g\big( (x,a,t,u),\Psi \big) \,:=\, \1\{x\in Q_W\}\, |f(\Delta_W)|,
\end{align*}
where $\Delta_W$ is regarded as a function of $\Psi$, yields
\begin{align*}
&\gamma\int_0^1 \int_{Q_W} \int_\AA \int_\MM \E{\big|f\big(\Delta_W^{(x,a,t,u)}\big)\big|} \Q \, (\d {u}) \; \Theta \, (\d {a})\; \d {x} \; \d t \\
& \qquad\qquad\qquad\qquad\qquad =\, \E{|f(\Delta_W)|\, \Psi(Q_W\times\AA\times [0,1]\times \MM) } \,<\, \infty.
\end{align*}
The finiteness follows from the Cauchy-Schwarz inequality.
Therefore, for almost all $(x,a,t,u)\in Q_W\times\AA\times [0,1]\times\MM$, the expectation $\E{\big|f\big(\Delta_W^{(x,a,t,u)}\big)\big|}$ is finite, and together with \eqref{Eq:bedEW} we obtain
\begin{align*}
\E{ D_{(x,a,t,u)}f(\Delta_W) \;\big|\; \Psi_t} \,&=\, \E{ f\big(\Delta_W^{(x,a,t,u)}\big) \;\big|\; \Psi_t} - \E{ f(\Delta_W) \;\big|\; \Psi_t} \\
&=\, \E{ f\big(\Delta_W^{(x,a,t,u)}\big) \;\big|\; \Psi_t} - \E{ f\big(\Delta_{W\setminus\{x\}}^{(x,a,t,u)}\big) \;\big|\; \Psi_t} \\
&=\, \E{ \Lambda_{(x,a,t,u)}f(\Delta_W) \;\big|\; \Psi_t}
\end{align*}
almost surely.
Since $\Lambda_{(x,a,t,u)}f(\Delta_W)$ vanishes whenever $x\notin W$, the statement follows.
\end{proof}

With the help of Lemma~\ref{Lem:exactVar} we can deduce a Poincaré type inequality (compare Theorem~5.2 in \cite{Last.RCM}).

\begin{Lem} \label{Lem:Poincare}
    Let $W\subseteq\R^d$ be measurable and bounded and $f$ a function defined for finite simplicial complexes with vertices in $\R^d\times\AA$, such that $\E{f(\Delta_W)^2}<\infty$.
    Then for all $t\in [0,1]$ it is true that
\begin{align*}
\V{f(\Delta_W)} \,\leq\, \gamma\int_{\R^d}\int_\AA\int_\MM \mathbb{E}\Big[ \left(\Lambda_{(x,a,t,u)}f(\Delta_W)\right)^2 \Big] \; \Q \, (\d {u}) \; \Theta \, (\d {a})\; \d {x}.
\end{align*}
\end{Lem}

\begin{proof}
Applying Jensen's inequality for conditional expectations to the integral representation of $\V{f(\Delta_W)}$ from Lemma~\ref{Lem:exactVar} yields the desired inequality by observing that the distribution of $\Lambda_{(x,a,t,u)}f(\Delta_W)$ does not depend on $t$.
\end{proof}

\section{Proof of the CLT and auxiliary results} \label{Sec:CLT}

The goal of this section is to prove the CLT (Theorem~\ref{Th:ZGWS3}).
The approach is based on the derivation of Theorem~2.4 in \cite{Can}, whose proof technique is, in turn, based on that of Theorem~3.1 in \cite{Trinh}.
In \cite{Can}, however, the difference operator used there is constructed slightly differently from the operator in (\ref{Def:DiffOp}), as the underlying Poisson process is first restricted to $W$, and then the simplicial complex is constructed from this restriction.
Theorem~2.4 in \cite{Can} is a CLT for weakly stabilizing functionals in the RCM as a random graph in the unmarked stationary case.
The proof is based on a general result for normal approximation of random variables, which can be found as Lemma~2.2 in \cite{Trinh} and will be stated here again for completeness.
For that, denote by $N_{\sigma^2}$, $\sigma^2 \in [0, \infty)$, a normally distributed random variable with mean 0 and variance $\sigma^2$.

\begin{Lem} [Lemma~2.2 in \cite{Trinh}] \label{Th:NormalApproxRV}
Let $(Y_n)_{n\in\N}$ and $(X_{n,L})_{n,L\in\N}$ be sequences of random variables with mean 0 satisfying the following two properties.
\begin{enumerate}
\item[$\mathrm{(i)}$] For all $L\in\N$ it is true that
\begin{align*}
\V{X_{n,L}} \,\rightarrow\, \sigma_L^2\,\in\, [0,\infty) \enskip\;\;\; \text{and} \;\;\;\enskip X_{n,L} \,\xrightarrow{d}\, N_{\sigma_L^2} \qquad \text{as } n \rightarrow \infty.
\end{align*}
\item[$\mathrm{(ii)}$] We have $\;\lim\limits_{L\rightarrow\infty} \limsup\limits_{n\rightarrow\infty}\, \V{ Y_n-X_{n,L} }=0$.
\end{enumerate}
Then the limit $\sigma^2:=\lim\limits_{L\rightarrow\infty} \sigma_L^2$ exists and
\begin{align*}
\V{Y_n} \,\rightarrow\, \sigma^2 \enskip\;\;\; \text{and} \;\;\;\enskip Y_n \,\xrightarrow{d}\, N_{\sigma^2} \qquad \text{as } n \rightarrow \infty
\end{align*}
holds.
\end{Lem}

Before proving the CLT, we first analyze the class of weakly
stabilizing functionals (Definition~\ref{Def:weaklyStabFunctionals}).
Recall that $V\sim\Theta$ and $U\sim\Q$ are independent random elements, and in particular independent of $\Psi$.
Throughout this section, we will consider the difference operator $\Lambda_{(\0,V,1,U)}f(\Delta_{W})$, defined in \eqref{Def:DiffOp}, where $W$ is a cube and $f$ a weakly stabilizing functional according to Definition~\ref{Def:weaklyStabFunctionals}.
From now on, we use the notation $|W|:=\lambda_d(W)$ for a measurable set $W\subseteq\R^d$.

\begin{Bem} \label{Kor:quadratischeIntegrierbarkeit}
Let $f$ be a weakly stabilizing functional satisfying the moment condition (\ref{MomentCondition}).
Then $\E{ f(\Delta_W)^2}<\infty$ holds for all $W\in\W$.
\end{Bem}

Remark~\ref{Kor:quadratischeIntegrierbarkeit} can be shown by approximating the function $f$ with bounded functions and applying the Poincaré-type inequality from Lemma~\ref{Lem:Poincare}.
A proof can be found in \cite{Pabst.Thesis} (Corollary~5.4).
The next proposition is analogous to the result in Lemma~2.7 in \cite{Can}, but we use a different argument.

\begin{Prop} \label{Prop:Lambda_t}
Let $f$ be a weakly stabilizing functional.
Then for each $t\in [0,1]$ there is a real-valued random variable $Z_t$ such that for all sequences of cubes $(W_n)_{n\in\N}$ with $W_n\rightarrow\R^d$ the convergence
\begin{align*}
\Lambda_{(\0,V,t,U)} \fWn \;\xrightarrow{\P}\; Z_t
\end{align*}
holds as $n\rightarrow\infty$.
\end{Prop}

\begin{proof}
Let $(W_n)_{n\in\N}$ be a sequence of cubes with $W_n\rightarrow\R^d$ and $t\in [0,1]$.
As a consequence of Lemma~5.6 in \cite{Kallenberg} we get
\begin{align*}
\E{\min\big\{|\Lambda_{(\0,V,1,U)}f(\Delta_{W_m})-\Lambda_{(\0,V,1,U)}f(\Delta_{W_n})|,1\big\}} \;\rightarrow\; 0 \qquad \text{ for } m,n\rightarrow\infty.
\end{align*}
Since for every choice of $n,m\in\N$ the equality in distribution
\begin{align*}
\left( \Lambda_{(\0,V,1,U)} \fWn, \Lambda_{(\0,V,1,U)} f\big(\Delta_{W_m}\big) \right) \;\stackrel{d}{=}\; \left( \Lambda_{(\0,V,t,U)} \fWn, \Lambda_{(\0,V,t,U)} f\big(\Delta_{W_m}\big) \right)
\end{align*}
holds, Lemma~5.6 in \cite{Kallenberg} implies that the sequence $(\Lambda_{(\0,V,t,U)}f(\Delta_{W_n}))_{n\in\N}$ converges in probability.
Since this is true for every such sequence of cubes, one can show in precisely the same manner as in the proof of Lemma~2.13 in \cite{Can}, that the limit variables according to different sequences of cubes have to be the same $\P$-almost surely.
\end{proof}

We proceed with the proof of Lemma~\ref{Lem:asymptoticVariance}.
The argument follows the same basic structure as the proof of Lemma~2.8 in \cite{Can}, though certain details differ, among other reasons, because the difference operator used in \cite{Can} is constructed somewhat differently.
For this reason, and to enhance clarity, we present the proof here, including details omitted in \cite{Can}.
We begin with an auxiliary result that provides an explicit representation of the asymptotic variance.

We decompose the probability space $\Omega = \Omega_1 \times \Omega_2$, where $(\Omega_1,\mathcal{A}_1,\P_1)$ is the probability space on which the Poisson process $\Psi$ is defined, and $(\Omega_2,\mathcal{A}_2,\P_2)$ is the probability space on which the marks $V, U$ are defined.
We denote the expectation with respect to $(\Omega_1,\mathcal{A}_1,\P_1)$ by $\tilde{\mathbb{E}}$.
Finally, let $\Psi_t$ be the restriction of $\Psi$ to $\R^d\times\AA\times [0,t)\times\MM$.

\begin{Lem} \label{Lem:VarRepresentation}
Let $f$ be a weakly stabilizing functional satisfying the moment condition~\eqref{MomentCondition}.
For every sequence of cubes $(W_n)_{n\in\N}$ with $W_n\rightarrow\R^d$, it is true that
\begin{align*}
\frac{\V{\fWn}}{|W_n|} \;\rightarrow\; \gamma\int_0^1 \E{ \Et{ Z_t \,|\, \Psi_t }^2 } \; \d t \,<\, \infty \qquad \text{as } n\rightarrow\infty,
\end{align*}
where the random variables $Z_t$, $t\in [0,1]$, are the random variables appearing in Proposition~\ref{Prop:Lambda_t}.
\end{Lem}

\begin{proof}
We fix $t\in [0,1]$ and a sequence of cubes $(W_n)_{n\in\N}$ with $W_n\rightarrow\R^d$.
The moment condition (\ref{MomentCondition}), which also applies to the operator $\Lambda_{(\0,V,t,U)}$, provides the uniform integrability of the sequence $(|\Lambda_{(\0,V,t,U)} f(\Delta_{W_n})|^2)_{n\in\N}$.
From Theorem~5.12 in \cite{Kallenberg}, Proposition~\ref{Prop:Lambda_t}, and the moment condition (\ref{MomentCondition}), we obtain
\begin{align} \label{L2_convergence}
\E{ \big| \Lambda_{(\0,V,t,U)} \fWn - Z_t \big|^2} \;\rightarrow\; 0 \qquad \text{as } n\rightarrow\infty
\end{align}
and also $\E{Z_t^2}<\infty$.
Jensen's inequality for conditional expectations, combined with the dominated convergence theorem, yields
\begin{align*}
&\int_0^1 \E{ \Et{ \Lambda_{(\0,V,t,U)} \fWn - Z_t \,\big|\, \Psi_t }^2 } \; \d t \\
&\qquad\qquad\qquad\qquad\qquad\leq\,  \int_0^1 \E{ \big| \Lambda_{(\0,V,t,U)} \fWn - Z_t \big|^2} \; \d t \;\rightarrow\; 0
\end{align*}
and therefore
\begin{align*}
h(W_n) \,:=\, \int_0^1 \E{ \Et{ \Lambda_{(\0,V,t,U)} \fWn \,\big|\, \Psi_t }^2 } \; \d t \;\rightarrow\; \int_0^1 \E{ \Et{ Z_t \,|\, \Psi_t }^2 } \; \d t.
\end{align*}
The integral on the right-hand side is finite and will be denoted by $c$ in the following.
Since this convergence holds for every sequence of cubes $(W_n)_{n\in\N}$ with $W_n\rightarrow\R^d$, we can find for every $\varepsilon>0$ a radius $r_\varepsilon>0$ such that
\begin{align} \label{EpsilonRadius}
\sup_{\substack{W \in \mathcal{W} \\ B(\0,r_\varepsilon)\subseteq W}}
\left| h(W) - c \right|
< \varepsilon.
\end{align}
Lemma~\ref{Lem:exactVar} provides for the functional $f(\Delta_{W_n})$
\begin{align*}
\V{f\big(\Delta_{W_n}\big)} \,&=\, \gamma \int_{W_n} \int_0^1 \E{ \Et{ \Lambda_{(x,V,t,U)}f\big(\Delta_{W_n}\big) \;\big|\; \Psi_t}^2 } \; \d t \; \d {x} \\
&=\, \gamma \int_{W_n} h(W_n-x) \; \d {x}.
\end{align*}
The last equality follows from
\begin{align*}
	\Et{ \Lambda_{(\0,V,t,U)} f\big(\Delta_{W_n-x}\big) \,\big|\, \Psi_t } \,\stackrel{d}{=}\, \Et{ \Lambda_{(x,V,t,U)} \fWn \,\big|\, \Psi_t }
\end{align*}
which is a consequence of the translation-invariant framework, where the distribution of $\Psi$, the connection functions, and the functional $f$ are invariant under translations with respect to the Euclidean component.
For $\varepsilon>0$, we choose $r_\varepsilon>0$ according to (\ref{EpsilonRadius}) and obtain
\begin{align} \label{ConvergenceVariance}
\frac{\V{\fWn}}{|W_n|} &= \frac{\gamma}{|W_n|} \int_{W_n} \1\big\{ d(x,\partial W_n)> r_\varepsilon \big\} h(W_n-x) \, \d x \nonumber \\
& \qquad + \frac{\gamma}{|W_n|} \int_{W_n} \1\big\{ d(x,\partial W_n)\leq r_\varepsilon \big\} h(W_n-x) \, \d x.
\end{align}
Jensen's inequality and the moment condition (\ref{MomentCondition}) yield, for $x\in W_n$,
\begin{align*}
h(W_n-x) \,&\leq\, \int_0^1 \E{ \big |\Lambda_{(\0,V,t,U)} f\big(\Delta_{W_n-x}\big)\big|^2 } \; \d t \\
&=\, \E{ \big |\Lambda_{(\0,V,1,U)} f\big(\Delta_{W_n-x}\big)\big|^2 } \,\leq\, \sup\limits_{\substack{W\in\W, \\ \0\in W}} \E{ \big|\Lambda_{(\0,V,1,U)} f(\Delta_W)\big|^{2} } \,<\, \infty,
\end{align*}
where we used that the distribution of $\Lambda_{(\0,V,t,U)}f\big(\Delta_{W_n-x}\big)$ does not depend on $t\in [0,1]$.
Hence, the second integral in (\ref{ConvergenceVariance}) converges to zero as $n\rightarrow\infty$ (compare (38) in \cite{Pabst.Euler} for more details).
On the other hand, using (\ref{EpsilonRadius}), it follows that
\begin{align*}
& \bigg| \, \frac{1}{|W_n|} \int_{W_n} \1\big\{ d(x,\partial W_n)> r_\varepsilon \big\} h(W_n-x) \, \d x - c \, \bigg| \\
& \qquad \leq\, \frac{1}{|W_n|} \int_{W_n} \1\big\{ d(x,\partial W_n)> r_\varepsilon \big\}\, \big|h(W_n-x)-c\big| \; \d x \\
& \qquad\qquad + \frac{c}{|W_n|}\, \big|\big\{x\in W_n \,|\, d(x,\partial W_n)\leq r_\varepsilon\big\}\big| \\
& \qquad \leq\, \varepsilon + \frac{c}{|W_n|}\, \big|\big\{x\in W_n \,|\, d(x,\partial W_n)\leq r_\varepsilon\big\}\big|.
\end{align*}
As the second term also asymptotically vanishes as $n\rightarrow\infty$ and $\varepsilon>0$ was arbitrary, the desired convergence follows.
\end{proof}

\begin{Lem} \label{Lem:AuxConvergence}
Assume that $f$ is a weakly stabilizing functional for which the moment condition~\eqref{MomentCondition} is fulfilled.
Then the following statements hold.
\begin{enumerate}
\item[$\mathrm{(i)}$] For every cube $W\in\W$, we have in $L^2$
\begin{align*}
	\Et{ \Lambda_{(\0,V,1,U)} f(\Delta_W) \,\big|\, \Psi_t } \;\rightarrow\; \Lambda_{(\0,V,1,U)} f(\Delta_W) \qquad\text{in } L^2\,\text{ as } t \to 1,
\end{align*}
\item[$\mathrm{(ii)}$] The following convergence holds in $L^2$
\begin{align*}
	\Et{ Z \,|\, \Psi_t } \;\rightarrow\; Z \qquad\text{ as } t \to 1.
\end{align*} 
\end{enumerate}
\end{Lem}

Before proceeding with the proof, we introduce the following notation.
We denote by $T(\eta)$ the simplicial complex constructed from a counting measure $\eta$ on $\R^d\times\AA\times [0,1]\times\MM$ according to Section~\ref{Sec:Model_Definition}, without formally introducing a space of simplicial complexes.
Recall that this construction is based on a decomposition $(Q_i)_{i\in\N_0}$ of $\R^d$ into cubes and $Q_0=[0,1)^d$ is the cube containing the origin $\0$.

\begin{proof}
We first prove (i).
Let $\xi_t$ be the restriction of the Poisson process $\Psi$ to $\R^d \times \AA \times [t,1] \times \MM$, and denote the expectation with respect to $\xi_t$ by $\mathbb{E}_t[\cdot]$.  
Furthermore, write
\begin{align*}
	\Lambda_{(\0,V,1,U)} f\big(T(\eta)_W\big) \,:=\, f\Big( T\big( \eta + \delta_{(\0,V,1,U)} \big)_W \Big) - f\Big( T\big( \eta + \delta_{(\0,V,1,U)} \big)_{W \setminus \{\0\}} \Big),
\end{align*}  
for a counting measure $\eta$ on $\R^d \times \AA \times [0,1] \times \MM$ and let $A_t$ be the event that $\xi_t$ has no points in $(W\cup Q_0)\times\AA\times [t,1]\times \MM$.
Observe that $\Lambda_{(\0,V,1,U)} f\big(T(\Psi)_W\big)=\Lambda_{(\0,V,1,U)} f\big(\Delta_W\big)$ and that $A_t$ is independent of $\Psi_t$.
On $A_t$, we have
\begin{align*}
	T\big( \Psi_t + \xi_t + \delta_{(\0,V,1,U)} \big)_W \,=\, T\big( \Psi_t + \delta_{(\0,V,1,U)} \big)_W \qquad \P\text{-almost surely}
\end{align*}
and hence
\begin{align} \label{AlmostSureEqual}
	\Lambda_{(\0,V,1,U)} f(\Delta_W) \,=\, \Lambda_{(\0,V,1,U)} f\big(T(\Psi_t)_W\big) \qquad \P\text{-almost surely}.
\end{align}
The conditional expectation appearing in (i) can be expressed as  
\begin{align*}
\P(A_t)\; \Lambda_{(\0,V,1,U)} f\big(T(\Psi_t)_W\big) \,+\, \mathbb{E}_t\left[ \Lambda_{(\0,V,1,U)} f\big(T(\Psi_t + \xi_t)_W\big)\, \1_{A_t^c} \right],
\end{align*}  
and thus, we obtain  
\begin{align}
&\Et{ \Lambda_{(\0,V,1,U)} f(\Delta_W) \,|\, \Psi_t } - \Lambda_{(\0,V,1,U)} f(T(\Psi_t)_W) \nonumber \\
&\hspace*{3pt}=\, (\P(A_t)-1)\; \Lambda_{(\0,V,1,U)} f(T(\Psi_t)_W) \,+\, \mathbb{E}_t\left[ \Lambda_{(\0,V,1,U)} f\big(T(\Psi_t + \xi_t)_W\big)\, \1_{A_t^c} \right]. \label{StochastischeKonvergenz_bedingte Erwartung_II}
\end{align}
Due to \eqref{AlmostSureEqual} and $\P(A_t)\rightarrow 1$ as $t\to1$, the convergence
\begin{align*}
\Lambda_{(\0,V,1,U)} f\big(T(\Psi_t)_W\big) \;\xrightarrow{\P}\; \Lambda_{(\0,V,1,U)} f\big(\Delta_W\big) \qquad \text{as } t \to 1
\end{align*}
holds and therefore, the first term in (\ref{StochastischeKonvergenz_bedingte Erwartung_II}) converges to 0 in probability.  
Applying Hölder's inequality (to the expectation with respect to $\xi_t$), we can also deduce that the second term in (\ref{StochastischeKonvergenz_bedingte Erwartung_II}) converges to 0 in probability.
Thus, the convergence in (i) holds in probability.
Together with the moment condition~\eqref{MomentCondition}, this implies convergence in $L^2$.

We now turn to (ii).
First, let a sequence of cubes $(W_n)_{n\in\N}$ with $W_n\rightarrow\R^d$ be given.  
Using Minkowski's inequality and the moment condition (\ref{MomentCondition}), we have  
\begin{align} \label{limsup_endlich}
\sup\limits_{n\in\N}\, \big\lVert\, \Et{ \Lambda_{(\0,V,1,U)} f(\Delta_{W_n}) \,|\, \Psi_t } - \Lambda_{(\0,V,1,U)} f(\Delta_{W_n}) \,\big\lVert_{L^2} \,<\,\infty.
\end{align}  
Additionally, Jensen's inequality for conditional expectations yields  
\begin{align*}
\big\lVert\, \Et{ \Lambda_{(\0,V,1,U)} f(\Delta_{W_n}) \,|\, \Psi_t } - \Et{ Z \,|\, \Psi_t } \,\big\lVert_{L^2} \,\leq\, \big\lVert\, \Lambda_{(\0,V,1,U)} f(\Delta_{W_n}) - Z \,\big\lVert_{L^2} \,\rightarrow\, 0
\end{align*}  
as $n \to \infty$.  
This gives us  
\begin{align*}
\big\lVert\, \Et{ Z \,|\, \Psi_t } - Z \,\big\lVert_{L^2} \,&\leq\, \limsup\limits_{n\rightarrow\infty}\, \big\lVert\, \Et{ Z \,|\, \Psi_t } - \Et{ \Lambda_{(\0,V,1,U)} f(\Delta_{W_n}) \,|\, \Psi_t }\,\big\lVert_{L^2} \\
&\enskip + \limsup\limits_{n\rightarrow\infty}\, \big\lVert\, \Et{ \Lambda_{(\0,V,1,U)} f(\Delta_{W_n}) \,|\, \Psi_t } - \Lambda_{(\0,V,1,U)} f(\Delta_{W_n})\,\big\lVert_{L^2} \\
&\enskip + \limsup\limits_{n\rightarrow\infty}\, \big\lVert\, \Lambda_{(\0,V,1,U)} f(\Delta_{W_n}) - Z\,\big\lVert_{L^2} \\
&=\, \limsup\limits_{n\rightarrow\infty}\, \big\lVert\, \Et{ \Lambda_{(\0,V,1,U)} f(\Delta_{W_n}) \,|\, \Psi_t } - \Lambda_{(\0,V,1,U)} f(\Delta_{W_n})\,\big\lVert_{L^2}.
\end{align*}  
Since the last limit superior is finite by (\ref{limsup_endlich}), for any $\varepsilon > 0$, there exists $N \in \N$ such that  
\begin{align*}
&\limsup\limits_{n\rightarrow\infty}\, \big\lVert\, \Et{ \Lambda_{(\0,V,1,U)} f(\Delta_{W_n}) \,|\, \Psi_t } - \Lambda_{(\0,V,1,U)} f(\Delta_{W_n})\,\big\lVert_{L^2} \\
&\qquad\qquad\leq\, \big\lVert\, \Et{ \Lambda_{(\0,V,1,U)} f(\Delta_{W_N}) \,|\, \Psi_t } - \Lambda_{(\0,V,1,U)} f(\Delta_{W_N})\,\big\lVert_{L^2} + \varepsilon.
\end{align*}  
From part~(i), it follows that  
\begin{align*}
\limsup\limits_{t\rightarrow 1}\, \big\lVert\, \Et{ Z \,|\, \Psi_t } - Z \,\big\lVert_{L^2} \,\leq\, \varepsilon.
\end{align*}  
Since $\varepsilon > 0$ was chosen arbitrarily, the assertion follows.
\end{proof}

With Lemmas~\ref{Lem:VarRepresentation} and~\ref{Lem:AuxConvergence} at hand, we now turn to the proof of Lemma~\ref{Lem:asymptoticVariance}.

\begin{proof}[Proof of Lemma~\ref{Lem:asymptoticVariance}]
The desired convergence of the variance follows directly from Lemma~\ref{Lem:VarRepresentation} by choosing
\begin{align*}
	\sigma^2 \,:=\, \gamma \int_0^1 \E{ \Et{ Z_t \mid \Psi_t }^2 } \, \d t.
\end{align*}
It remains to show that $\sigma^2 > 0$ whenever $\P(Z \neq 0) > 0$.
To this end, we demonstrate the continuity of the integrand in the definition of $\sigma^2$ at $t = 1$, that is
\begin{align} \label{ContinuityIntegrand}
\E{ \Et{ Z_t \,|\, \Psi_t }^2 } \;\rightarrow\; \E{ Z^2 } \qquad \text{as } t \to 1.
\end{align}
Since $\E{Z^2}>0$ whenever $\P(Z \neq 0) > 0$, the convergence in \eqref{ContinuityIntegrand} implies the positivity of $\sigma^2$.
To conclude the proof, it remains to prove \eqref{ContinuityIntegrand}.
For this purpose, consider the event $B_t$ that $\Psi$ has no points in $Q_0\times\AA \times[t,1]\times\MM$.
By the construction of the simplicial complex, on the event $B_t$, we have  
\begin{align*}
T\big( \Psi + \delta_{(\0,V,t,U)} \big)_W \, = \, T\big( \Psi + \delta_{(\0,V,1,U)} \big)_W \qquad \P\text{-almost surely}
\end{align*}
for every cube $W\in\W$ and therefore, in particular,  
\begin{align*}
\Lambda_{(\0,V,t,U)} f\big(\Delta_{W}\big) \, = \, \Lambda_{(\0,V,1,U)} f\big(\Delta_{W}\big) \qquad \P\text{-almost surely}.
\end{align*}  
Since the probability of the event $B_t$ converges to 1 as $t \to 1$, we obtain
\begin{align} \label{StochastischeKonvergenz_tnach1}
\Lambda_{(\0,V,t,U)} f\big(\Delta_{W}\big) \;\xrightarrow{\P}\; \Lambda_{(\0,V,1,U)} f\big(\Delta_{W}\big)\qquad \text{as } t \to 1
\end{align}
for every cube $W\in\W$.
The moment condition in fact yields convergence in $L^2$.
Now consider a sequence of cubes $(W_n)_{n\in\N}$ with $W_n\to\R^d$.
The triangle inequality implies for each $n\in\N$ that  
\begin{align*}
&\Vert Z_t - Z \Vert_{L^2} \,\leq\, \Vert Z_t - \Lambda_{(\0,V,t,U)} \fWn \Vert_{L^2} \\
& \qquad + \Vert \Lambda_{(\0,V,t,U)} \fWn - \Lambda_{(\0,V,1,U)} \fWn \Vert_{L^2} + \Vert \Lambda_{(\0,V,1,U)} \fWn - Z \Vert_{L^2}.
\end{align*}  
The first and last term on the right-hand side vanish, due to (\ref{L2_convergence}), as $n \to \infty$.
Using that the convergence in (\ref{StochastischeKonvergenz_tnach1}) also holds in $L^2$, we thus obtain the convergence $Z_t\to Z$ as $t\to 1$. 
Together with Jensen's inequality for conditional expectations, this yields 
\begin{align*}
\E{ \Et{ Z_t - Z \,|\, \Psi_t }^2 } \;\rightarrow\; 0 \qquad \text{as } t \to 1.
\end{align*} 
Together with part~(ii) of Lemma~\ref{Lem:AuxConvergence}, the convergence \eqref{ContinuityIntegrand} follows and the proof is completed.
\end{proof}

Using Lemma~\ref{Lem:asymptoticVariance} and Lemma~\ref{Th:NormalApproxRV}, we prove the CLT for weakly stabilizing functionals stated in Theorem~\ref{Th:ZGWS3}.

\begin{proof}[Proof of Theorem~\ref{Th:ZGWS3}]
Fix a sequence $(W_n)_{n\in\N}$ with $W_n\rightarrow\R^d$.
The distributions of $f(\Delta_{W+t})$ and $f(\Delta_W)$ are the same for all cubes $W\in\W$ and $t\in\R^d$ due to the translation-invariant framework ($f$, the connection functions, and the distribution of $\Psi$ are translation invariant with respect to the Euclidean component).
Therefore, we may assume that all cubes of the sequence $(W_n)_{n\in\N}$ are centered at the origin.
Hence, there exists a sequence $(m_n)_{n\in\N}$ of positive numbers with $m_n\rightarrow\infty$ such that
\begin{align*}
	W_n \,=\, \Big[ -\tfrac{\sqrt[d]{m_n}}{2}, \tfrac{\sqrt[d]{m_n}}{2} \Big)^d,
\end{align*} 
that is, $m_n$ is the volume of $W_n$.
Now fix $L \in \mathbb{N}$ and decompose the space into cubes of side length $L$, i.e., consider the decomposition $([0, L)^d + L \cdot z_i)_{i \in \mathbb{N}_0}$ with the previously chosen enumeration $z_0, z_1, \dots$ of the lattice $\mathbb{Z}^d$ with $z_0 = \0$.
Let $l_n\in\N_0$ and $(G_i)_{i \in \mathbb{N}}$ be an enumeration of the cubes in this decomposition, such that $G_i \subseteq W_n$ holds precisely for $i\leq l_n$.
Since $L$ is a natural number, each cube $G_i$ is a union of cubes from the decomposition $(Q_i)_{i\in\N_0}$.
We define a double sequence, which will play the role of the double sequence in Lemma~\ref{Th:NormalApproxRV}, by
\begin{align*}
X_{n,L} \,&:=\, \frac{1}{\sqrt{m_n}} \sum_{i=1}^{l_n} \fD{G_i} - \mathbb{E} \left[ \fD{G_i} \right].
\end{align*}
Recall that, by the construction of $\Delta$ introduced in Section~\ref{Sec:Model_Definition}, the subcomplex $\Delta_{G_i}$ is almost surely completely determined by the process $\Psi_{G_i}$ and the processes $\Psi_{G_i}$, $i\in\N$, are independent.
Therefore, from the latter representation, it follows immediately that $X_{n,L}$ is a sum of i.i.d. random variables.
Moreover, Remark~\ref{Kor:quadratischeIntegrierbarkeit} yields that these random variables have finite variance.
It follows that
\begin{align*}
\V{X_{n,L}} = \frac{l_n}{m_n} \, \V{\fD{G_1}} \;\rightarrow\; \frac{1}{L^d} \V{\fD{G_1}} \,=:\, \sigma_L^2 \,\in\, [0,\infty)
\end{align*}
as $n\rightarrow\infty$.
The convergence $\frac{l_n}{m_n} \rightarrow \frac{1}{L^d}$ follows from the second convergence in (9.15) in \cite{Schneider}.
Note that by Lemma~\ref{Lem:asymptoticVariance}, we have $\sigma^2 = \lim_{L \rightarrow \infty} \sigma_L^2$.
Thus, the sequence of random variables $(X_{n,L})$ satisfies a CLT as $n \rightarrow \infty$ (cf. e.g., Theorem~27.1 in \cite{Billingsley}), i.e., it holds
\begin{align*}
X_{n,L} \;\xrightarrow{d}\; N_{\sigma_L^2} \qquad \text{as } n \rightarrow \infty.
\end{align*}
Thus, the double sequence $(X_{n,L})_{n,L \in \mathbb{N}}$ satisfies the first condition of Lemma~\ref{Th:NormalApproxRV}.
To prove the second condition, let $Y_n := \frac{1}{\sqrt{m_n}}\left( f\big(\Delta_{W_n}\big) - \mathbb{E} \left[ f\big(\Delta_{W_n}\big) \right] \right)$.
Then, by the Poincaré inequality (Lemma~\ref{Lem:Poincare}), we have
\begin{align}
&\V{Y_n -X_{n,L}} \nonumber \\
&\qquad\leq\, \frac{\gamma}{m_n} \int_{W_n} \E{\bigg( \Lambda_{(x,V,1,U)} \fWn - \sum_{i=1}^{l_n} \1\{ x\in G_i \} \; \Lambda_{(x,V,1,U)} f\big(\Delta_{G_i}\big) \bigg)^2 } \; \d x \nonumber \\
&\qquad=\, \frac{\gamma}{m_n} \int_{W_n\setminus (\cup_{i=1}^{l_n} G_i)} \E{\left( \Lambda_{(x,V,1,U)} \fWn \right)^2 } \; \d x \nonumber \\
&\qquad\qquad + \frac{\gamma}{m_n} \sum_{i=1}^{l_n} \int_{G_i} \E{\left( \Lambda_{(x,V,1,U)} \fWn - \Lambda_{(x,V,1,U)} f\big(\Delta_{G_i}\big) \right)^2 } \; \d x. \label{Varlimsup}
\end{align}
The expectation in the first summand of (\ref{Varlimsup}) is uniformly bounded in $n$ due to the moment condition (\ref{MomentCondition}).
From the first equation in (9.15) in \cite{Schneider}, it follows that the first summand asymptotically vanishes as $n \rightarrow \infty$.
The convergence (\ref{L2_convergence}) from the proof of Lemma~\ref{Lem:asymptoticVariance} provides for $t = 1$
\begin{align*}
\E{ \big| \Lambda_{(\0,V,1,U)} \fWn - Z \big|^2} \;\rightarrow\; 0.
\end{align*}
This convergence holds for any sequence of cubes $(W_n)_{n\in\N}$ with $W_n\rightarrow\R^d$.
Therefore, for any $\varepsilon > 0$, there exists a radius $r_\varepsilon > 0$ such that it is true that
\begin{align} \label{EpsilonRadius2}
\E{ \big| \Lambda_{(\0,V,1,U)} f(\Delta_W) - \Lambda_{(\0,V,1,U)} f(\Delta_V) \big|^2} \,<\, \varepsilon
\end{align}
for all $W, V \in \W$ with $B(\0, r_\varepsilon) \subset W \cap V$.
For a fixed $\varepsilon > 0$, we define for a cube $W \in \W$
\begin{align*}
\mathrm{int}_\varepsilon (W) \,&:=\, \big\{ x \in W \,\|\, d(x,\partial W) \geq r_\varepsilon \big\}, \\
\partial_\varepsilon (W) \,&:=\, \big\{ x \in W \,\|\, d(x,\partial W) < r_\varepsilon \big\} \,=\, W \setminus \mathrm{int}_\varepsilon (W).
\end{align*}
For $L > 2r_\varepsilon$, we have $|\partial_\varepsilon G_i| = L^d - (L - 2r_\varepsilon)^d$, and for sufficiently large $L$, we get $|\partial_\varepsilon G_i| \,\leq\, 2r_\varepsilon d L^{d-1}$.
Due to the distributional equality
\begin{align*}
\big(\Lambda_{(\0,V,1,U)} f(\Delta_W), \Lambda_{(\0,V,1,U)} f(\Delta_V)\big) \stackrel{d}{=} \big(\Lambda_{(x,V,1,U)} f(\Delta_{W+x}), \Lambda_{(x,V,1,U)} f(\Delta_{V+x})\big)
\end{align*}
for all cubes $W,V \in \W$ and $x \in \R^d$, it follows from (\ref{EpsilonRadius2})
\begin{align*}
\E{ \big| \Lambda_{(x,V,1,U)} \fD{W_n} - \Lambda_{(x,V,1,U)} \fD{G_i} \big|^2} \,<\, \varepsilon \qquad \text{for all } \; x \in \mathrm{int}_\varepsilon(G_i), \enskip i \leq l_n.
\end{align*}
The moment condition (\ref{MomentCondition}) ensures the existence of a constant $\rho > 0$ such that
\begin{align*}
\E{ \big| \Lambda_{(x,V,1,U)} f(\Delta_W) - \Lambda_{(x,V,1,U)} f(\Delta_V) \big|^2 } \,\leq\, \rho \qquad \text{for all } \; W,V \in \W, \; x \in \R^d.
\end{align*}
Thus, for the second summand in (\ref{Varlimsup}), using the previous three estimates for sufficiently large $L$, we obtain
\begin{align*}
\frac{\gamma}{m_n} \sum_{i=1}^{l_n} \int_{G_i} &\E{\left( \Lambda_{(x,V,1,U)} \fWn - \Lambda_{(x,V,1,U)} f\big(\Delta_{G_i}\big) \right)^2 } \; \d x \\
&\leq\, \frac{\gamma}{m_n} \sum_{i=1}^{l_n} \bigg( \int_{\mathrm{int}_\varepsilon(G_i)} \E{\left( \Lambda_{(x,V,1,U)} \fWn - \Lambda_{(x,V,1,U)} f\big(\Delta_{G_i}\big) \right)^2 } \; \d x \\
&\qquad + \int_{\partial_\varepsilon(G_i)} \E{\left( \Lambda_{(x,V,1,U)} \fWn - \Lambda_{(x,V,1,U)} f\big(\Delta_{G_i}\big) \right)^2 } \; \d x \bigg) \\
&\leq\, \gamma \frac{l_n}{m_n} \Big( \varepsilon (L - 2r_\varepsilon)^d + 2 \rho r_\varepsilon d L^{d-1} \Big) \\
&\leq\, \gamma\varepsilon + 2 \gamma\rho r_\varepsilon d L^{-1},
\end{align*}
where we used the trivial estimates $l_n L^d \leq m_n$ and $(L - 2r_\varepsilon)^d \leq L^d$ for the last inequality.
We conclude that
\begin{align*}
\limsup\limits_{L\rightarrow\infty} \limsup\limits_{n\rightarrow\infty}\, \V{ Y_n - X_{n,L}} \,\leq\, \gamma\varepsilon.
\end{align*}
Since $\varepsilon$ was arbitrary, condition (ii) in Lemma~\ref{Th:NormalApproxRV} follows, and thus the desired distributional convergence holds.
Lemma~\ref{Lem:asymptoticVariance} provides the additional result on the positivity of the asymptotic variance $\sigma^2$.
\end{proof}

Finally, from Theorem~\ref{Th:ZGWS3}, we also directly obtain a multivariate CLT for weakly stabilizing functionals.

\begin{Kor} [Multivariate CLT for weakly stabilizing functionals]
Let $m \in \N$ and $g_1, \dots, g_m$ be weakly stabilizing functionals that satisfy the moment condition (\ref{MomentCondition}).
According to Proposition~\ref{Prop:Lambda_t}, there exist random variables $Z^{(i)}_t$, $i \in \{1, \dots, m\}$, $t \in [0, 1]$, such that
\begin{align*}
\Lambda_{(\0, V, t, U)} g_i\big(\Delta_{W_n}\big) \;\xrightarrow{\P}\; Z_t^{(i)} \qquad \text{as } n \to \infty
\end{align*}
for all sequences of cubes $(W_n)_{n \in \N}$ with $W_n \to \R^d$.
Then, for all $i,j \in \{1,\dots,m\}$ and any sequence of cubes $(W_n)_{n \in \N}$ with $W_n \to \R^d$,
\begin{align} \label{asymptoticCov_schwachstabilisierend}
\lim\limits_{n \to \infty} \frac{\Cov{g_i\big(\Delta_{W_n}\big)}{g_j\big(\Delta_{W_n}\big)}}{|W_n|} \,=\, \gamma\int_0^1 \mathbb{E}\Big[ \mathbb{E}\big[ Z_t^{(i)} \,|\, \Psi_t \big] \, \mathbb{E}\big[ Z_t^{(j)} \,|\, \Psi_t \big] \Big] \; \d t \,=:\, \sigma_{i,j}
\end{align}
and the covariance matrix $\Sigma := (\sigma_{i,j})_{1\leq i,j\leq m}$ is positive semi-definite.
Let $N_\Sigma$ be an $m$-dimensional centered normally distributed random vector with covariance matrix $\Sigma$.
Then
\begin{align*}
\frac{1}{\sqrt{|W_n|}} \left( g_1\big(\Delta_{W_n}\big) - \E{g_1\big(\Delta_{W_n}\big)}, \dots, g_m\big(\Delta_{W_n}\big) - \E{g_m\big(\Delta_{W_n}\big)} \right) \;\xrightarrow{d}\; N_\Sigma
\end{align*}
as $n\to\infty$.
\end{Kor}

\begin{proof}
The claim follows by applying Theorem~\ref{Th:ZGWS3}, since a linear combination of weakly stabilizing functionals that satisfy the moment condition (\ref{MomentCondition}) is itself weakly stabilizing and fulfills the moment condition.
For further details, refer to the proof of Corollary~2.12 in \cite{Can}.
\end{proof}

\section{Betti numbers and other weakly stabilizing functionals} \label{Sec:weaklyStabFct}

After deriving a CLT for an abstract class of functionals in the previous section, we will now focus on specific example functionals in this section.
We will show that many functionals from the literature on random simplicial complexes or random graphs, such as Betti numbers or the Euler characteristic, satisfy the conditions of Theorem~\ref{Th:ZGWS3} under certain integrability conditions.
To demonstrate that these functionals are weakly stabilizing, we will use a characterization of weak stabilization that simplifies the definition in two aspects.
Using this characterization, we can even establish almost sure convergence of the difference operator $\Lambda_{(\0,V,1,U)}$.
Throughout this section, when referring to the added vertex, we write $\0$ instead of $(\0,V)$.
However, before we turn to example functionals, we will examine the edge degree $\deg_1(\0,\Delta^\0)$ of the origin in $\Delta^\0:=\Delta^{(\0,V,1,U)}$.
This will be the key quantity for verifying the moment condition for the functionals, as it will allow us to control the difference operator.

It is well known that the edge degree of a deterministic vertex added to the RCM, even on more general spaces, follows a Poisson distribution.
More precisely, adding a deterministic vertex $(\0,a)$ for some fixed $a\in\AA$ would lead to a Poisson distribution of the edge degree of the added vertex with parameter
\begin{align} \label{Mixing_Parameter}
    \pi(a) \,:=\, \gamma \int_{\R^d} \int_\AA \varphi_1\big( (\0,a),(y,b) \big) \; \Theta (\d b) \; \d y.
\end{align}
Therefore, by choosing a random mark $V$ of the added vertex, which is independent of the remaining part of the model, we gain a mixed Poisson distribution with mixing parameter $\pi(V)$.
A mixed Poisson distribution can be considered as the result of a two-stage random experiment, where in the first sub-experiment the random parameter is drawn, and in the second sub-experiment a Poisson-distributed random variable is generated, with the outcome of the first sub-experiment serving as the parameter.
Mixed Poisson distributions also appear in generalized random graphs (see, for example, Section~6.3 in \cite{vanderHofstad}).
The precise definition and some stochastic properties of a mixed Poisson distribution can be found in \cite{Kuba}, where the crucial statement for our work is Proposition~1.
This states, among other things, that the $s$-th moment of a mixed Poisson distribution exists if the $s$-th moment of its mixing parameter exists, i.e.
\begin{align} \label{Momentimplication}
    \E{\pi(V)^s} \,<\,\infty \enskip\Longrightarrow\enskip \E{\deg_1(\0,\Delta^\0)^s} \,<\,\infty, \qquad s\in\N.
\end{align}
We will use this statement to verify the moment condition \eqref{MomentCondition} for the considered functionals under suitable moment assumptions on the mixing parameter $\pi(V)$.
Note that in the unmarked case, the mixing parameter is constant and therefore all its moments exist, as long as this value is assumed to be finite, which leads to the typical integrability condition.
Next, we will formulate the aforementioned characterization of weak stabilization.

\begin{Prop} [Characterization of weak stabilization] \label{Prop:CharacterizationWS}
For a translation-invariant functional $f$ defined on finite simplicial complexes with vertices in $\R^d\times\AA$ the following statements are equivalent.
\begin{enumerate}
\item[$\mathrm{(i)}$] $f$ is weakly stabilizing.
\item[$\mathrm{(ii)}$] The operator $\Lambda_{(\0,V,1,U)}f(\Delta_{W_n})$ converges for all increasing sequences of cubes $(W_n)_{n\in\N}$ with $W_n\rightarrow\R^d$ in probability as $n\rightarrow\infty$.
\end{enumerate}
\end{Prop}

\begin{proof}
The argumentation works in exactly the same way like the argumentation in the proof of Lemma~2.13 in \cite{Can}.
\end{proof}

If we want to show that a certain functional $f$ is weakly stabilizing, this is simplified in two ways by Proposition~\ref{Prop:CharacterizationWS}.
On the one hand, we only need to consider increasing sequences $(W_n)_{n\in\N}$ of cubes, which in turn provides us with an increasing sequence of simplicial complexes $(\Delta^\0_{W_n})_{n\in\N}$.
On the other hand, for each such sequence, it is only necessary to prove the existence of the limit of $\Lambda_{(\0,V,1,U)}f(\Delta_{W_n})$, and not that these limits agree for all sequences.

\subsection{Betti numbers}

In this subsection, we aim to prove, using Proposition~\ref{Prop:CharacterizationWS}, that the Betti numbers defined in Subsection~\ref{Sec:Definition_Betti} satisfy the conditions of Theorem~\ref{Th:ZGWS3} whenever an integrability assumption for the mixing parameter $\pi(V)$ is satisfied.

\begin{Prop} \label{Betti:schwachstabilisierend}
Let $p\in\N_0$ such that $\E{\pi(V)^{3(p+1)}}<\infty$ holds.
Then the $p$-th Betti number $\beta_p$ is a weakly stabilizing functional satisfying the moment condition (\ref{MomentCondition}).
\end{Prop}

\begin{proof}
First of all, $\beta_p$ is clearly translation invariant since it depends only on the combinatorial structure of a simplicial complex.
To show that the convergence (\ref{Convergence_weaklyStab}) holds, we fix an increasing sequence $(W_n)_{n\in\N}$ of cubes with $W_n\rightarrow\R^d$ and a realization $\omega\in\Omega$ such that $\deg_1(\0,\Delta^\0)$ is finite, which holds with probability one.
Then $(\Delta^\0_{W_n})_{n\in\N}$ is an increasing sequence of finite simplicial complexes and Lemma~3.8 in \cite{Can} shows that the limit
\begin{align*}
Z \,:=\, \lim\limits_{n\rightarrow\infty} \Lambda_{(\0,V,1,U)}\beta_p\big(\Delta_{W_n}\big) \,=\, \lim\limits_{n\rightarrow\infty} \beta_p(\Delta^\0_{W_n})-\beta_p(\Delta^\0_{W_n\setminus\{\0\}})
\end{align*}
exists.
Hence, Proposition~\ref{Prop:CharacterizationWS} implies that $\beta_p$ is weakly stabilizing.
To verify the moment condition (\ref{MomentCondition}), we apply Lemma~\ref{Lem:Betti_properties} (iv) with $K=\Delta^\0_W,L=\Delta^\0_{W\setminus\{\0\}}$ and obtain
\begin{align*}
\sup\limits_{\substack{W\in\W, \\ \0\in W}} \E{ \big|\Lambda_{(\0,V,1,U)} \beta_p(\Delta_W)\big|^{3} } \,&\leq\, \sup\limits_{\substack{W\in\W, \\ \0\in W}} \E{ \big| \deg_p(\0,\Delta^\0_W) + \deg_{p+1}(\0,\Delta^\0_W) \big|^{3} } \\
&\leq\, \E{ \big| \deg_p(\0,\Delta^\0) + \deg_{p+1}(\0,\Delta^\0) \big|^{3} } \\
&\leq\, \mathbb{E} \Bigg[ \bigg| \binom {\deg_1(\0,\Delta^\0)} p + \binom {\deg_1(\0,\Delta^\0)} {p+1} \bigg|^{3} \Bigg] \\
&\leq\, \mathbb{E} \Bigg[ \bigg| \deg_1(\0,\Delta^\0)^p +\deg_1(\0,\Delta^\0)^{p+1} \bigg|^{3} \Bigg] \,<\,\infty.
\end{align*}
For the third inequality, we used the fact that every $p$-simplex that contains the origin can be identified with a $p$-element set of neighbors of the origin.
Due to (\ref{Momentimplication}) and the assumptions, the upper bound is finite.
\end{proof}

For $\alpha=1$ and $p=1$, Proposition~\ref{Betti:schwachstabilisierend} yields a CLT for the first Betti number of the random graph, which has not been established so far, even in the unmarked case.
Finally, we aim to provide a law of large numbers for the Betti numbers.
For this, we require the integrability condition
\begin{align} \label{Integrabilität:esssup}
\esssup\limits_{a\in\AA} \; \int_{\R^d}\int_\AA \varphi_1\big((\0,a),(y,b)\big) \; \Theta(\d b) \; \d y \,<\, \infty.
\end{align}
Note that (\ref{Integrabilität:esssup}) implies that the mixing parameter $\pi(V)$ is bounded $\P$-almost surely.

\begin{Kor}
Suppose that condition (\ref{Integrabilität:esssup}) holds.
Then there exists a constant $b_p\geq 0$ for $p\in\N_0$, such that for all sequences of cubes $(W_n)_{n\in\N}$ with $W_n\rightarrow\R^d$, we have
\begin{align*}
\frac{\beta_p\big(\Delta_{W_n}\big)}{|W_n|} \;\xrightarrow{\P}\; b_p.
\end{align*}
\end{Kor}

\begin{proof}
We fix a sequence of cubes $(W_n)_{n\in\N}$ with $W_n\rightarrow\R^d$.
First, the convergence of $\mathbb{E}\big[\beta_p\big(\Delta_{W_n}\big)\big]|W_n|^{-1}$ follows in the same way as in the proof of Theorem~3.6 in \cite{Can}, and we set $b_p$ as the limit of this sequence.
In the proof of Theorem~3.6 in \cite{Can}, properties (iii),(iv) from Lemma~\ref{Lem:Betti_properties} of the Betti numbers were used, as well as the existence of the limit 
\begin{align*}
    \lim_{n\rightarrow\infty} \frac{\E{ f_i(\Delta_{W_n}) }}{|W_n|}
\end{align*}
for each $i\in\{0,\dots,\alpha\}$, where $f_i(K)$ denotes the number of $i$-simplices of a simplicial complex $K$.
The latter is shown for the present situation in Section~7 of \cite{Pabst.Euler}, where condition (\ref{Integrabilität:esssup}) is used.
Since, according to Lemma~\ref{Lem:asymptoticVariance} and Proposition~\ref{Betti:schwachstabilisierend}, the limit of \( \mathbb{V}\big(\beta_p\big(\Delta_{W_n}\big)\big)|W_n|^{-1} \) also exists, the convergence in probability follows directly from Chebyshev's inequality.
The assumption from Proposition~\ref{Betti:schwachstabilisierend} is satisfied, as (\ref{Integrabilität:esssup}) implies the almost sure boundedness of $\pi(V)$.
\end{proof}

\subsection{Euler characteristic} \label{Sec:Euler}

Having established that the Betti numbers fall into the class of functionals to which Theorem~\ref{Th:ZGWS3} is applicable, in this section we will turn to other example functionals.
Besides the Betti numbers, the Euler characteristic is one of the most popular topological quantities when studying random simplicial complexes (see, for example, \cite{Candela,Thomas}).
Recall that $f_i(K)$ denotes the number of $i$-simplices of a simplicial complex $K$.
For a finite simplicial complex $K$, the Euler characteristic is defined by
\begin{align*}
    \chi(K) \,:=\, \sum_{i=0}^{\dim(K)} (-1)^i f_i(K).
\end{align*}
Like the Betti numbers, the Euler characteristic can be defined on much more general topological spaces.
The CLT in \cite{Hug} for the Boolean model can be applied, for example, to the Euler characteristic, which is found there as the zeroth intrinsic volume.

\begin{Kor} \label{Kor:Euler_weakstabilizing}
Assume $\E{\pi(V)^{3\alpha}}<\infty$.
Then the Euler characteristic $\chi$ is weakly stabilizing and the moment condition (\ref{MomentCondition}) holds for $\chi$.
\end{Kor}

\begin{proof}
To apply Proposition~\ref{Prop:CharacterizationWS}, we fix an increasing sequence of cubes $(W_n)_{n\in\N}$ with $W_n\rightarrow\R^d$.
The difference operator from (\ref{Def:DiffOp}) for the Euler characteristic satisfies
\begin{align} \label{Euler_DiffOp}
\Lambda_{(\0,V,1,U)}\chi\big(\Delta_{W_n}\big) \,&=\, \sum_{i=0}^\alpha (-1)^i \deg_i(\0,\Delta_{W_n}^\0).
\end{align}
Since $\E{\pi(V)}<\infty$ implies $\deg_1(\0,\Delta^\0)<\infty$, and therefore $\deg_i(\0,\Delta^\0)<\infty$ for all $i\leq\alpha$, almost surely, the right-hand side of (\ref{Euler_DiffOp}) converges almost surely as $n\rightarrow\infty$.
Hence, according to Proposition~\ref{Prop:CharacterizationWS}, $\chi$ is weakly stabilizing.
Furthermore, (\ref{Euler_DiffOp}) yields
\begin{align*}
\sup\limits_{\substack{W\in\W, \\ \0\in W}} \E{ \big|\Lambda_{(\0,V,1,U)} \chi(\Delta_W)\big|^{3} } \,&\leq\, \sum_{i,j,m=0}^\alpha \; \E{ \deg_i(\0,\Delta^\0)\deg_j(\0,\Delta^\0)\deg_m(\0,\Delta^\0) } \\
&\leq\, \sum_{i,j,m=0}^\alpha \; \E{ \deg_1(\0,\Delta^\0)^{i+j+m} }.
\end{align*}
Here, for the second inequality, we used $\deg_k(\0,\Delta^\0)\leq \deg_1(\0,\Delta^\0)^k$ (compare the proof of Proposition~\ref{Betti:schwachstabilisierend}).
Due to the condition $\E{\pi(V)^{3\alpha}}<\infty$ and (\ref{Momentimplication}), all moments appearing in the last expression exist, and the upper bound is therefore finite.
\end{proof}

Since the Euler characteristic can also be expressed as a linear combination of Betti numbers according to the Euler-Poincaré formula
\begin{align*}
    \chi(K) \,=\, \sum_{p=0}^{\dim(K)} (-1)^p \, \beta_p(K),
\end{align*}
it could be directly inferred from Proposition~\ref{Betti:schwachstabilisierend} that the Euler characteristic is weakly stabilizing and satisfies the moment condition (\ref{MomentCondition}).
However, this would lead to the stronger condition $\E{\pi(V)^{3(\alpha+1)}}<\infty$.
This is because, in the proof of Proposition~\ref{Betti:schwachstabilisierend}, the $p$-th Betti number is estimated from above using $f_p$ and $f_{p+1}$.
For $p=\alpha$, however, the term $f_{p+1}$ does not appear, which leads to a weaker moment condition.

As mentioned in the introduction, even when using Corollary~\ref{Kor:Euler_weakstabilizing}, the CLT cannot be applied to the Euler characteristic of the \v Cech or Vietoris–Rips complex, since their Euler characteristics do not coincide with those of their $\alpha$-skeletons.
As both complexes are uniquely determined by their vertex sets, we may define their Euler characteristics as functions of finite point sets, which in turn yield functionals on finite simplicial complexes through their vertex set.
We briefly introduce both complexes and omit the mark space, as it is not needed here and would only complicate the notation.
The \v{C}ech complex with parameter $r>0$ on a vertex set $X$ contains a simplex  consisting of the points $x_0,\dots,x_n\in X$ if and only if
\begin{align*}
   \bigcap_{i=0}^n B(x_i,r) \,\neq\, \emptyset,
\end{align*}
where $B(x,r)$ denotes the ball of radius $r$ centered at $x$.
For the Vietoris–Rips complex, the corresponding condition is
\begin{align*}
	\Vert x_i-x_j\Vert \,\leq\, r \qquad\text{for all } i\neq j.
\end{align*}
For a finite set $X\subset\mathbb{R}^d$, we define $H_{\mathrm{C}}^r(X)$ and $H_{\mathrm{VR}}^r(X)$ as the Euler characteristics of the \v Cech and Vietoris–Rips complexes with parameter $r>0$ on the vertex set $X$, respectively.
In this way, we obtain functionals on finite simplicial complexes through
\begin{align*}
	\chic(K) \,:=\, H_{\mathrm{C}}^r(\Vrt(K)), \qquad \chivr(K) \,:=\, H_{\mathrm{VR}}^r(\Vrt(K)),
\end{align*}
where $\Vrt(K)$ denotes the vertex set of $K$.
It follows immediately that $\chic$ and $\chivr$ are translation invariant.
Since these definitions only depend on the vertex set of the simplicial complex, the choice of $\alpha$ and the connection functions are irrelevant.

\begin{Kor}
The functionals $\chic$ and $\chivr$ are weakly stabilizing and satisfy the moment condition~(\ref{MomentCondition}).
\end{Kor}

\begin{proof}
Fix an increasing sequence of cubes $(W_n)_{n\in\N}$ with $W_n\rightarrow\R^d$ and observe that
\begin{align} \label{DiffOp_Euler_Cech}
	\Lambda_{(\0,V,1,U)}\chic\big(\Delta_{W_n}\big) \,=\, \sum_{i=0}^\infty (-1)^i \deg_i^\mathrm{C}(\0,\Delta_{W_n}^\0),
\end{align}
where $\deg_i^\mathrm{C}(\0,\Delta_{W_n}^\0)$ denotes the $i$-simplex degree of $\0$ in the \v Cech complex with parameter $r>0$ on the vertex set $\Vrt(\Delta_{W_n})\cup\{\0\}$.
Since in the \v Cech complex the vertex $\0$ can only be connected to points at distance at most $2r$, the operator in \eqref{DiffOp_Euler_Cech} becomes constant as $n \to \infty$ once $B(\0,2r) \subset W_n$.
Therefore, $\chic$ is weakly stabilizing, and the same argument shows that $\chivr$ is weakly stabilizing.
To verify the moment condition, let $N$ be the number of vertices of $\Delta$ in $B(\0,2r)$.
Thus, we obtain for $\chic$
\begin{align*}
	\sup\limits_{\substack{W\in\W, \\ \0\in W}} \E{ \big|\Lambda_{(\0,V,1,U)} \chic(\Delta_W)\big|^{3} } \,&\leq\, \E{ \Big|\sum_{i=0}^\infty \deg_i^\mathrm{C}(\0,\Delta) \Big|^3 } \,\leq\, \E{ \Big|\sum_{i=0}^\infty \binom N i \Big|^3 } \\
&\leq\, \E{ 8^{N} } \,=\, \exp\left( 7\gamma|B(\0,2r)| \right).
\end{align*}
The last equality can be computed directly using that $N$ follows a Poisson distribution with parameter $\gamma|B(\0,2r)|$.
The same bound applies to the difference operator of $\chivr$, and therefore both functionals satisfy the moment condition \eqref{MomentCondition}.
\end{proof}

\subsection{Subcomplex counts and further example functionals}

The Euler characteristic considered in the previous section is a linear combination of the functionals $f_i$ (number of $i$-simplices), i.e., quantities that count specific substructures, in this case simplices of a certain dimension.
As further example functionals, we want to consider functionals that count certain (possibly even more complex) substructures within simplicial complexes.

For this purpose, let $L$ be a connected, finite simplicial complex, and let $g_L$ be the functional that assigns to a simplicial complex the number of induced subcomplexes isomorphic to $L$.
An induced subcomplex of a simplicial complex $K$ is the restriction of this complex to simplices with vertices in a given subset $I$ of the vertex set of $K$, denoted by $K[I]$.
Writing $\Vrt(K)$ for the vertex set of a simplicial complex $K$, the functional $g_L$ can be written as
\begin{align*}
	g_L(K) \,:=\, \sum_{I\subseteq\Vrt(K)} \1\{ K[I]\cong L \}.
\end{align*}
Since $K[I]\cong L$ can only occur if $|I|$ equals the number of vertices of $L$, the sum could be restricted to subsets $I$ of this size.
However, as this is not necessary and would make the formulas more cumbersome, we will always sum over all subsets in the following.
We recover the functional $f_i$ as a special case when we choose $L$ as a simplicial complex consisting only of an $i$-simplex and all its subsimplices.
Furthermore, we denote by $h_L(K)$ the number of connected components of a simplicial complex $K$ isomorphic to $L$.

Functionals of this type can often be found in the literature.
For example, counts of components of a certain form (i.e., isomorphic to a fixed chosen graph) were already considered in \cite{Last.RCM}, while the CLT in \cite{Can} is applied to isomorphic subgraph counts, among other things.
In \cite{Penrose.IRG}, functionals counting induced connected subgraphs or connected components of a given size are studied.

\begin{Prop} \label{Prop:schwachStab}
Let $L$ be a connected simplicial complex on $r+1$ vertices for some $r\in\N_0$.
\begin{enumerate}
\item If the integrability condition (\ref{Integrabilität:esssup}) is satisfied, then $g_L$ is a weakly stabilizing functional that fulfills the moment condition (\ref{MomentCondition}).
\item In the case that $\E{\pi(V)^3}<\infty$, $h_L$ is a weakly stabilizing functional that meets the moment condition (\ref{MomentCondition}).
\end{enumerate}
\end{Prop}

\begin{proof}
For the proof, let us fix an increasing sequence of cubes $(W_n)_{n\in\N}$ with $W_n\rightarrow\R^d$.
Furthermore, we write $\Delta^\0[\0,I]$ for the restriction of $\Delta^\0$ to the set $I\cup\{\0\}$, i.e., to simplices with vertices in this set.
Considering the functional $g_L$, we first note that in the case $r=0$, the functional coincides with the number of vertices $f_0$.
In this case, the difference operator is identically equal to 1 once $\0\in W_n$, and the statement is therefore obvious.
In the following, let $r\in\N$.
The operator $\Lambda_{(\0,V,1,U)}g_L\big(\Delta_{W_n}\big)$ counts all induced subcomplexes of $\Delta_{W_n}^\0$ isomorphic to $L$ that contain the origin and is thus non-negative.
Concretely, the operator satisfies
\begin{align*}
\Lambda_{(\0,V,1,U)}g_L\big(\Delta_{W_n}\big) \,&=\, \sum_{I\subseteq \Vrt(\Delta_{W_n})} \1\big\{ \Delta_{W_n}^\0 [ \0,I ]\cong L \big\} \\ 
& \xrightarrow{a.s.} \enskip \sum_{I\subseteq \Vrt(\Delta)} \1\big\{ \Delta^\0 [ \0,I ]\cong L \big\} \,=:\, Z_L\qquad\text{as } n\rightarrow\infty.
\end{align*}
An application of Mecke’s formula (Theorem~4.4 in \cite{Last.Lectures}) shows that the expected value of $Z_L$ is equal to
\begin{align*}
\frac{\gamma^r}{r!} \int_{\left( \R^d \right)^{r}} \int_{\AA^{r+1}} \P\left( \Delta\big((\0,a_0),(x_1,a_1),\dots,(x_r,a_r)\big) \cong L \right) \; \Theta^{r+1} (\d {a}) \; \d {x}.
\end{align*}
where $\d a = \d {(a_0,\dots,a_r)}$ and $\d x = \d {(x_1,\dots,x_r)}$.
We can bound the integrand by $(r+1)!f_L((\0,a_0),(x_1,a_1),\dots,(x_r,a_r))$, where $f_L$ denotes the simplex function of $L$ according to Definition~4.1 in \cite{Pabst.Euler}.
Therefore, Corollary~7.2 in \cite{Pabst.Euler} implies that the expectation of $Z_L$ is finite, where condition (\ref{Integrabilität:esssup}) is applied.
Hence, $Z_L$ is almost everywhere finite and Proposition~\ref{Prop:CharacterizationWS} implies that $g_L$ is weakly stabilizing.
To prove the validity of the moment condition for the functional $g_L$, we show $\E{Z_L^3}<\infty$.
Initially, we have
\begin{align*}
    Z_L^3 \,=\, \sum_{I,J,M\subseteq \Vrt(\Delta)} \1\big\{ \Delta^\0 [ \0,I ], \Delta^\0 [ \0,J ], \Delta^\0 [ \0,M ]\cong L \big\}.
\end{align*}
Since $L$ is connected, we can bound the indicator function appearing here by the indicator function of the event that $\Delta^\0$ restricted to $I\cup J\cup M\cup\{\0\}$ is connected.
By distinguishing cases based on the cardinality of this set, we can find constants $C(s)$ such that
\begin{align*}
    Z_L^3 \,\leq\, \sum_{s=r}^{3r} \sum_{I\subseteq \Vrt(\Delta), |I|=s} C(s)\; \1\big\{ \Delta^\0[\0,I] \text{ is connected}\big\}.
\end{align*}
Clearly a simplicial complex on $s+1$ vertices is connected if and only if it is isomorphic to a connected simplicial complex with vertex set $\{1,\dots,s+1\}$ and there are only finitely many of them.
Denoting the set of such simplicial complexes by $B_s$, we obtain
\begin{align*}
    Z_L^3 \,\leq\, \sum_{s=r}^{3r} \sum_{I\subseteq \Vrt(\Delta), |I|=s}\, \sum_{K\in B_s} \;C(s)\; \1\big\{ \Delta^\0[\0,I]\cong K\big\} \,=\, \sum_{s=r}^{3r} \sum_{K\in B_s} \;C(s)\; Z_K.
\end{align*}
Following the same argument used previously for $Z_L$, we find that the expected values of the random variables $Z_K$ are also finite, proving statement (i).

To describe the difference operator of the functional $h_L$ for the second part of the proposition, let $Z(v,K)$ denote the connected component of a simplicial complex $K$ that contains the vertex $v\in K$.
Moreover, for a subset $I\subseteq\Vrt(K)$, we denote by $K[I]\stackrel{cc}{\cong} L$ the event that the restriction of $K$ to $I$ is a connected component of $K$ isomorphic to $L$.
Finally, we write $I\leftrightarrow v$ if there is an edge between a vertex of $I$ and $v$.
With this notation, we derive
\begin{align*}
\Lambda_{(\0,V,1,U)}h_L\big(\Delta_{W_n}\big) \,&=\, \1\big\{ Z\big(\0,\Delta_{W_n}^\0\big)\cong L \big\} - \sum_{I\subseteq \Vrt(\Delta_{W_n})} \1\big\{ \Delta_{W_n} [ I ]\stackrel{cc}{\cong} L, I\leftrightarrow \0 \big\} \\
&\xrightarrow{a.s.} \enskip \1\big\{ Z\big(\0,\Delta^\0\big)\cong L \big\} - \sum_{I\subseteq\Vrt(\Delta)} \1\big\{ \Delta[ I ]\stackrel{cc}{\cong} L, I\leftrightarrow \0 \big\}
\end{align*}
as $n\rightarrow\infty$.
To understand the almost sure finiteness of this limit variable, we observe that the sum involved is bounded above by $\deg_1(\0,\Delta^\0)$.
Thus, together with the assumption and (\ref{Momentimplication}), we also have
\begin{align*}
&\E{|\Lambda_{(\0,V,1,U)}h_L(\Delta_{W})|^3} \,\leq\, \E{ (1+\deg_1(\0,\Delta^\0))^3 } \,<\, \infty.
\end{align*}
\end{proof}

As the final example functional, we consider, for $m\in\N,l\in\N_0$, the number of vertices with $m$-simplex degree equal to $l$, denoted by $d_m^l$.
For a finite simplicial complex $K$, define
\begin{align*}
    d_m^l(K) \,:=\, \sum_{v\in\Vrt(K)} \1\big\{ \deg_m(v,K)=l \big\}
\end{align*}
with the definition~(\ref{SimplexDegree}) of a simplex degree.
This functional for $m=1$ appears in the literature for random graphs and is examined, for example, in Section~8.4 of \cite{vanderHofstad} for the Preferential Attachment Model.
In the RCM, vertex degree distributions are studied in \cite{Iyer.Degrees}, for example.

\begin{Kor}
Let $m\in\N$, $l\in\N_0$, and assume that $\E{\pi(V)^3}<\infty$.
Then $d_m^l$ is a weakly stabilizing functional that satisfies the moment condition (\ref{MomentCondition}).
\end{Kor}

\begin{proof}
We again fix an increasing sequence of cubes $(W_n)_{n\in\N}$ with $W_n\rightarrow\R^d$ to apply Proposition~\ref{Prop:CharacterizationWS}.
The operator $\Lambda_{(\0,V,1,U)}d_m^l\big(\Delta_{W_n}\big)$ can be written as
\begin{align*}
&\1\big\{ \deg_m(\0,\Delta_{W_n}^\0)=l \big\} \\
&\qquad+ \sum_{v\in\Vrt(\Delta_{W_n})} \1\big\{ \deg_m\big(v,\Delta_{W_n}^\0\big)=l \big\} - \1\big\{ \deg_m\big(v,\Delta_{W_n\setminus\{\0\}}^\0\big)=l \big\} \\
& \xrightarrow{a.s.} \enskip \1\big\{ \deg_m(\0,\Delta^\0)=l \big\} \\
& \qquad + \sum_{v\in\Vrt(\Delta)} \1\big\{ \deg_m(v,\Delta^\0)=l \big\} - \1\big\{ \deg_m\big(v,\Delta_{\R^d\setminus\{\0\}}^\0\big)=l \big\}.
\end{align*}
The limit variable is almost surely finite, as the absolute value of the sum involved can be bounded above by $\deg_1(\0,\Delta^\0)$.
With this reasoning, we also have $|\Lambda_{(\0,V,1,U)}d_m^l(\Delta_{W})| \leq 1 + \deg_1(\0,\Delta^\0)$, from which the desired statement follows together with the assumption and (\ref{Momentimplication}).
\end{proof}

\noindent Finally, we point out that, although only isomorphism-invariant example functionals were considered in this section, Theorem~\ref{Th:ZGWS3} also permits more general functionals.

\section{On the positivity of the asymptotic variance} \label{Sec:Variance}

For the functionals introduced in the previous section, Corollary~\ref{Kor:Positivity} gives a criterion for the strict positivity of the asymptotic variance in the CLT (Theorem~\ref{Th:ZGWS3}).
We now turn to its proof, which is based on the final part of Theorem~\ref{Th:ZGWS3}.

\begin{proof}[Proof of Corollary~\ref{Kor:Positivity}]
We will prove the result for the case in which the second statement holds, as the first statement implies this.
This follows because the function in the integral from (ii) is almost everywhere strictly positive in this case.
Now, let $K$ be as specified in statement (ii) and $A\subseteq\Omega$ the event that the connected component containing the origin in $\Delta^\0$ is isomorphic to $K$.
Using Mecke’s formula (Theorem~4.4 in \cite{Last.Lectures}), we find that the expected number of connected components in $\Delta^\0$ isomorphic to $K$ that contain the origin, and consequently $\P(A)$, is positive whenever the integral from (ii) is positive.

Now let $\omega\in A$ and $W\in\W$ be a sufficiently large cube such that all the vertices of the connected component of the origin in $\Delta^\0$ are contained in $W$.
Since this connected component is isomorphic to $K$, it is in particular finite, and we can always find such a cube.
By restricting the connected component of the origin to all vertices except the origin, we obtain a complex that is isomorphic to an induced subcomplex $L$ of $K$ on $r$ vertices.
Due to the last required property of $K$ and the weak additivity of $f$, it follows that
\begin{align*}
|\Lambda_{(\0,V,1,U)} f(\Delta_{W})| \,=\, |f(K) - f(L)| \,>\, 0.
\end{align*}
Therefore, the random variable $Z$ associated with $f$ (from Definition~\ref{Def:weaklyStabFunctionals}) is not zero on $A$, and in particular it follows that $\P(Z\neq 0)>0$, which, according to Theorem~\ref{Th:ZGWS3}, enforces the positivity of the asymptotic variance $\sigma^2$.
\end{proof}

For the application of Corollary~\ref{Kor:Positivity} to a functional $f$ with the required properties, the first question that arises is whether there exists a connected simplicial complex $K$ with $|f(K)-f(L)|>0$ for all induced subcomplexes with one fewer vertex, for which the conditions (i),(ii) could be verified in the case $\dim(K)\leq\alpha$.
If this is the case, then the asymptotic variance from Theorem~\ref{Th:ZGWS3} is always positive according to Corollary~\ref{Kor:Positivity}, provided that the connection functions take values only in $(0,1)$.
All functionals from Section~\ref{Sec:weaklyStabFct} are isomorphism-invariant and weakly additive.
We want to provide examples for the choice of the complex $K$ from Corollary~\ref{Kor:Positivity}, starting with the Betti numbers.
For the zeroth Betti number, a simplicial complex consisting of a single vertex can always be chosen.
A minimal example (in terms of the number of vertices) for a simplicial complex with a $p$-dimensional hole for $p\in\N$, i.e., with a positive $p$-th Betti number, is given by the family of all proper subsets of a set with $p+2$ elements, for example by
\begin{align} \label{Boundary_p+1-Simplex}
    K_p \,=\, \big\{ \sigma\subseteq \{1,\dots,p+2\} : 0<|\sigma|< p+2 \big\}. 
\end{align}
Since every simplicial complex with a positive $p$-th Betti number has at least $p+2$ vertices, $\beta_p(L)=0$ holds for all induced subcomplexes $L$ of $K_p$ on $p+1$ vertices.
In case $\varphi_{p+1}\equiv 1$, however, the two conditions from Corollary~\ref{Kor:Positivity} are not fulfilled for this simplicial complex.
In this situation, we construct a simplicial complex $L_p$ for which Condition (ii) from Corollary~\ref{Kor:Positivity} can be verified.
On the vertex set $\{1,\dots,p+3\}$, we consider the $p$-simplices
\begin{align*}
    \sigma_j \,:=\, \{1,\dots,p+2\}\setminus\{j\}, \qquad \rho_j \,:=\, \{1,\dots,p+1,p+3\}\setminus\{j\}
\end{align*}
for $j\in\{1,\dots,p+1\}$.
Let $L_p$ be the simplicial complex consisting of all simplices $\sigma_j,\rho_j$ for $j\in\{1,\dots,p+1\}$ and all their subsimplices, i.e.
\begin{align*}
    L_p \,=\, \big\{ \tau\subseteq \{1,\dots,p+3\} : \exists j\in\{1,\dots,p+1\} \text{ with } \emptyset\neq\tau\subseteq \sigma_j \text{ or } \emptyset\neq\tau\subseteq \rho_j \big\}.
\end{align*}
For the $p$-dimensional simplicial complex $L_p$, the properties $\beta_p(L_p)$=1 and $\beta_p(M)=0$ for all induced subcomplexes $M$ of $L_p$ on $p+2$ vertices can be verified.
Additionally, it does not contain the boundary of a $(p+1)$-simplex, meaning all its subsimplices without the $(p+1)$-simplex itself, which allows it to be applied when $\varphi_{p+1} \equiv 1$.
Another possible choice for this situation is the simplicial complex $O_p$ from Definition~3.3 in \cite{Kahle.RandomGeometric}.
Since $O_p$, however, contains $2p+2$ vertices and thus, in case $p \geq 2$, more than the simplicial complex $L_p$, the complex $L_p$ may be better suited for verifying Condition (ii) from Corollary~\ref{Kor:Positivity}.
In the case $p=1$, both complexes are isomorphic.

Considering the Euler characteristic $\chi$, the complex $K_1$ from (\ref{Boundary_p+1-Simplex}) can be used since $\chi(K_1)=0$ and $\chi(M)=1$ holds for all induced subcomplexes $M$ of $K_1$ on 2 vertices.
For the functionals $g_L,h_L$, the simplicial complex $L$ satisfies the desired property.
To find a suitable simplicial complex for $d_m^l$ (the number of vertices with $m$-simplex degree $l$), we ask whether there exists a connected simplicial complex $K$ with $\deg_m(v,K)=l$ for all vertices $v$ of $K$.
Such a simplicial complex would satisfy the desired property.
The existence of such a simplicial complex with dimension $m$ is equivalent to the existence of a 1-$(v,m+1,l)$-design (in the sense of Definition~12.1 in \cite{Jukna}) for some $v\in\N$.
Such a design exists according to Theorem~12.3 in \cite{Jukna} whenever we can find a $b\in\N$ with
\begin{align*}
    m+1 \,<\, v, \qquad b\,\leq\, \binom {v}{m+1}, \qquad b(m+1)=vl.
\end{align*}
This scenario occurs with $v=k(m+1)$ and $b=kl$, provided $k\geq 2$ is sufficiently large.
If the associated simplicial complex is not connected, each connected component satisfies the desired property.

\section{Application to the Boolean model} \label{Sec:BooleanModel}

\begin{figure}[t!]
  \captionsetup{ labelfont = {bf}, format = plain }
  \centering
  \subfloat[][]{\includegraphics[width=0.48\linewidth]{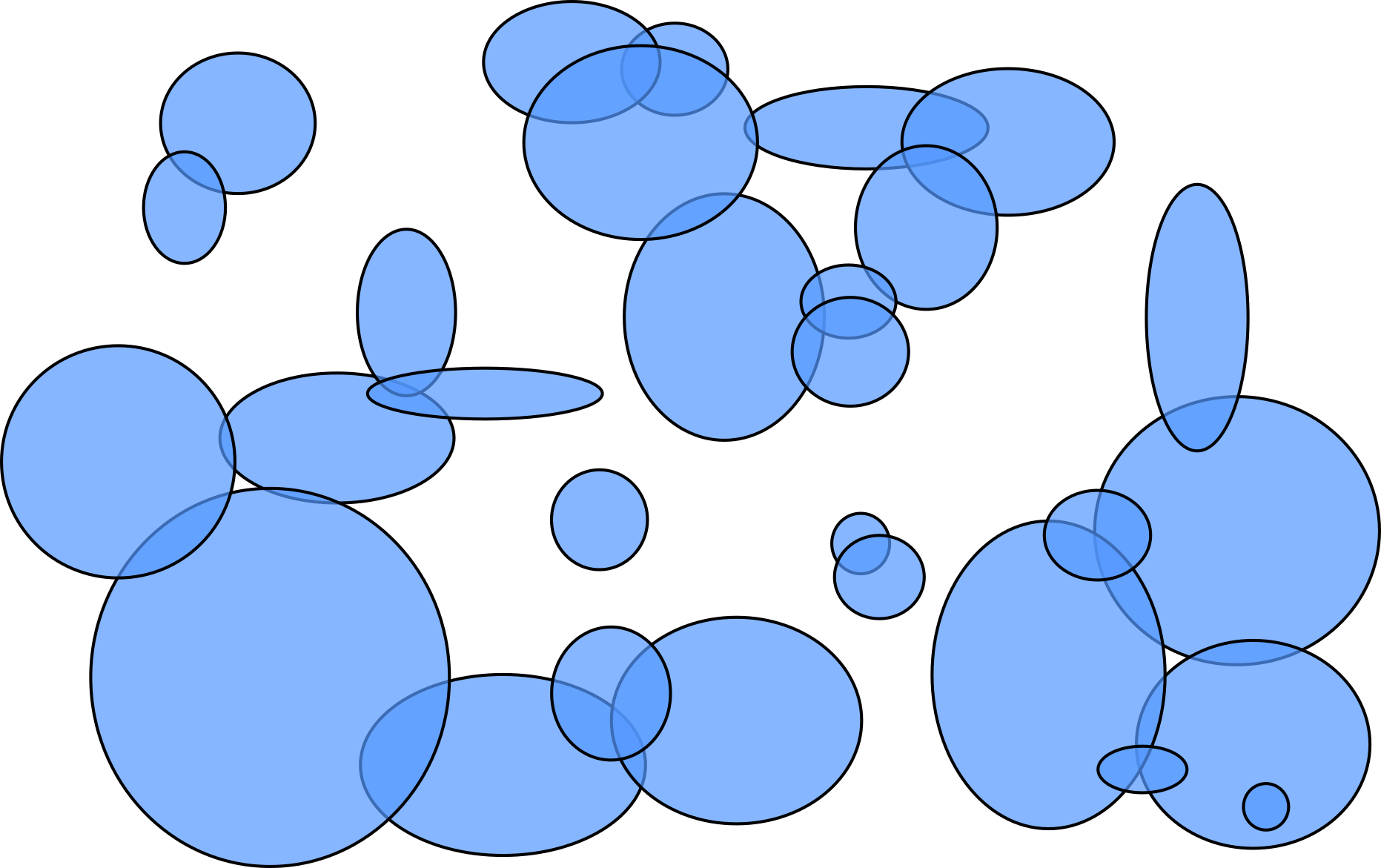}}
  \quad
  \subfloat[][]{\includegraphics[width=0.48\linewidth]{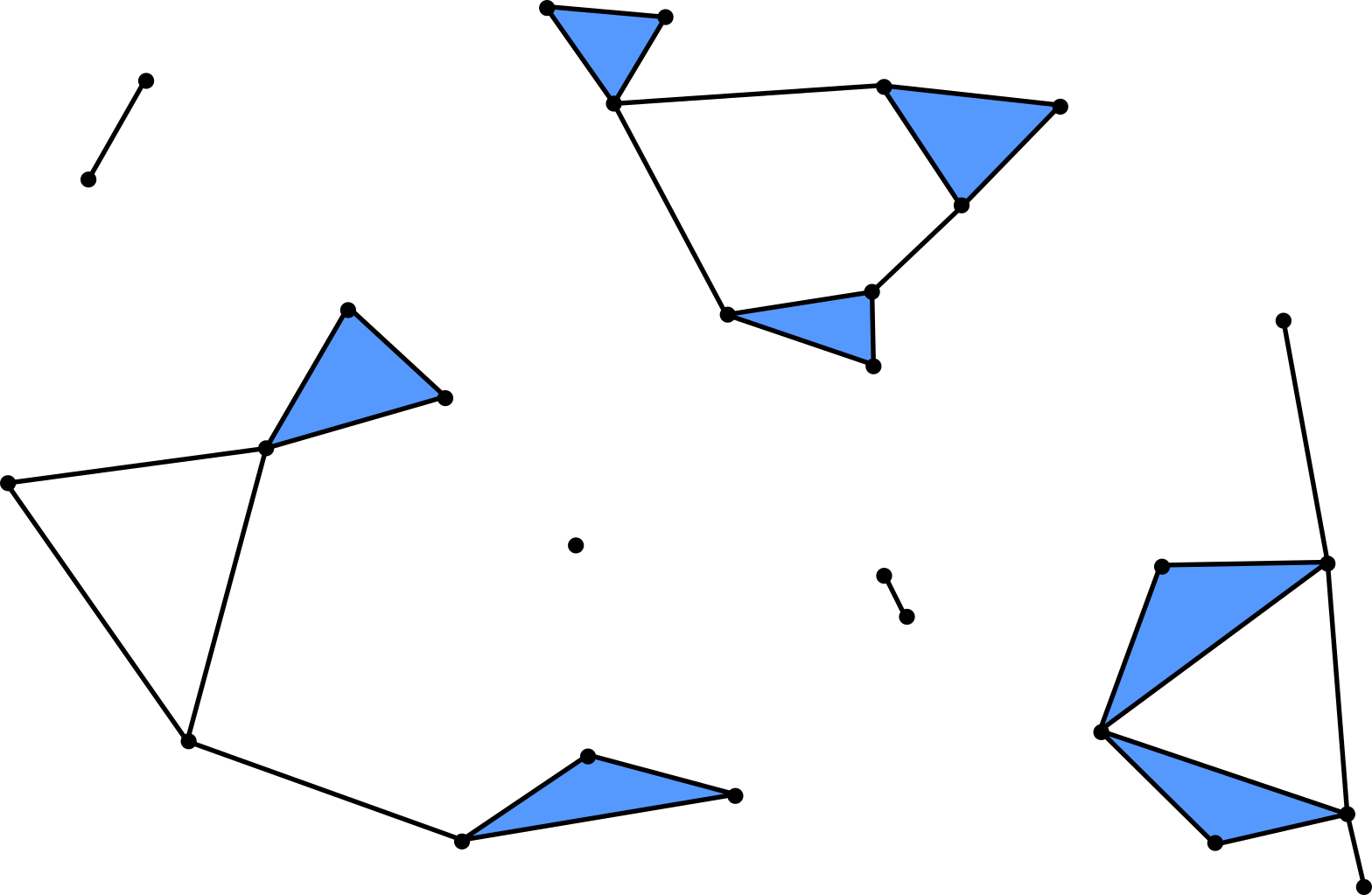}}
  \caption{(a) A union of convex sets and (b) the corresponding homotopy-equivalent nerve}
  \label{fig:BooleanNerve}
\end{figure}

Finally, in this section, we want to unravel the connection to the Boolean model that was already hinted at in the introduction, thereby obtaining a CLT for Betti numbers in the Boolean model.
To this end, we choose as the marking space the set $\AA=\mathcal{K}^d$ of all non-empty, compact and convex subsets of $\R^d$, equipped with the $\sigma$-field induced by the Hausdorff metric.
According to the Theorems~A.19 and~A.26 in \cite{Last.Lectures}, this is a Borel space.
As before, let $\Theta$ be a probability measure on $\mathcal{K}^d$. 
Additionally, we choose the connection functions 
\begin{align} \label{BoolM_connectionfcts}
    \varphi_j\big((x_0,K_0),\dots,(x_j,K_j)\big) = \1\Big\{ \bigcap_{i=0}^j (x_i+K_i)\neq\emptyset \Big\}, \qquad j\in \{1,\dots,\alpha\},
\end{align}
which are obviously translation invariant in the sense of (\ref{translation_invariant}).
Since these connection functions only take values in $\{0,1\}$, we can forgo the additional marking space $\MM$ and thus construct the random simplicial complex $\Delta$ directly from the Poisson process $\Phi$ on $\R^d\times\mathcal{K}^d$.
To establish the connection to the Boolean model, we first want to introduce it.
In stochastic geometry, the Boolean model (in $\R^d$) corresponding to the Poisson process $\Phi$ with intensity measure $\gamma\lambda_d\otimes\Theta$ refers to the random closed set 
\begin{align*}
    Z \,:=\, \bigcup_{(x,K)\in\Phi} (x+K)
\end{align*}
and the individual union sets are called the grains of the model.
The probability measure $\Theta$ is called the grain distribution and $\gamma$ the intensity of the Boolean model.
For a bounded, measurable set $W\subset\R^d$, define
\begin{align*}
    Z^W \,:=\, \bigcup_{(x,K)\in\Phi_W} (x+K),
\end{align*}
which is a finite union of convex sets.
Be aware of the difference compared to the restriction of $Z$ to $W$.
A variant of the Nerve Theorem (p.~76 in \cite{Edelsbrunner}) shows that the set $Z^W$ is homotopy equivalent to the nerve of the individual union sets.
The nerve of a finite collection $A$ of nonempty sets is the simplicial complex defined by
\begin{align*}
    N(A) \,:=\, \big\{ \emptyset\neq \sigma\subseteq A : \cap_{B\in\sigma} B\neq\emptyset \big\}.
\end{align*}
Figure~\ref{fig:BooleanNerve} shows a union of 29 convex sets and the corresponding nerve $N$ with $\beta_0(N)=6$, $\beta_1(N)=3$ and $\beta_p(N)=0$ for all $p\geq 2$.
Note that $\Delta_W$ is the $\alpha$-skeleton of the nerve of the grains $x+K$ with $x\in W$.
Consequently, it follows that for a functional $f$ defined on arbitrary topological spaces, which is invariant under homotopy equivalence and depends only on their $\alpha$-skeletons for simplicial complexes, the equality $f(Z^W)=f(\Delta_W)$ holds almost surely.
This is true for the $p$-th Betti number $\beta_p$ whenever $p\leq\alpha-1$ due to Lemma~\ref{Lem:Betti_properties} (ii).
Therefore, for $p\leq\alpha-1$, applying Theorem~\ref{Th:ZGWS3} to $\beta_p$ yields a CLT for $\beta_p(Z^W)$.
Since $\alpha\in\N$ can be chosen arbitrarily, the condition $p\leq\alpha-1$ imposes no restriction at all.

Finally, we take a look at the integrability condition of Proposition~\ref{Betti:schwachstabilisierend} in the Boolean model.
For the $k$-th moment of the mixing parameter $\pi(V)$ from (\ref{Mixing_Parameter}), $V\sim\Theta$, we obtain by applying Jensen's inequality
\begin{align}
    \E{\pi(V)^k} \,&=\, \int_{\mathcal{K}^d} \left( \gamma \int_{\R^d} \int_{\mathcal{K}^d} \1\{ K\cap (L+y)\neq\emptyset \} \; \Theta (\d L) \; \d y \right)^k \; \Theta (\d K) \nonumber \\
    &=\, \int_{\mathcal{K}^d} \left( \gamma \int_{\mathcal{K}^d} |K-L| \; \Theta (\d L) \right)^k \; \Theta (\d K) \nonumber \\
    &\leq\, \gamma^k \int_{\mathcal{K}^d} \int_{\mathcal{K}^d} |K-L|^k \; \Theta (\d L) \; \Theta (\d K), \label{Integrability_BoolM}
\end{align}
where $A-B:=\{ a-b\mid a\in A, b\in B\}$ denotes the Minkowski difference of two sets $A,B\subseteq\R^d$.
Proposition~\ref{Betti:schwachstabilisierend} states that Theorem~\ref{Th:ZGWS3} is applicable for the $p$-th Betti number, whenever the integral in (\ref{Integrability_BoolM}) is finite for $k=3(p+1)$.
We summarize the discussion of this section in the following corollary.

\begin{Kor}
Let $p\in\N_0$ and consider a Boolean model $Z$ with grain distribution $\Theta$ and intensity $\gamma$ satisfying
\begin{align} \label{BoolM_BettiInt}
	\int_{\mathcal{K}^d} \int_{\mathcal{K}^d} |K-L|^{3(p+1)} \; \Theta (\d L) \; \Theta (\d K) \,<\, \infty.
\end{align}
Then there exists $\sigma^2\geq 0$ such that for any sequence of cubes $(W_n)_{n\in\N}$ with $W_n\rightarrow\R^d$ it is true that
\begin{align*}
\frac{\beta_p(Z^{W_n}) - \E{ \beta_p(Z^{W_n}) }}{\sqrt{|W_n|}} \;\xrightarrow{d}\; N_{\sigma^2} \qquad \text{as } n \to \infty.
\end{align*}
\end{Kor}

\begin{proof}
Choose $\alpha=p+1$, the mark space $\AA=\mathcal{K}^d$ and the connection functions from \eqref{BoolM_connectionfcts}.
The condition \eqref{BoolM_BettiInt} ensures that $\E{\pi(V)^{3(p+1)}}<\infty$ as demonstrated in \eqref{Integrability_BoolM}.
Hence, Proposition~\ref{Betti:schwachstabilisierend} shows that the $p$-th Betti number $\beta_p$ satisfies the requirements of the CLT.
Theorem~\ref{Th:ZGWS3} therefore yields 
\begin{align*}
\frac{\beta_p\big(\Delta_{W_n}\big) - \E{ \beta_p\big(\Delta_{W_n}\big) }}{\sqrt{|W_n|}} \;\xrightarrow{d}\; N_{\sigma^2} \qquad \text{as } n \to \infty,
\end{align*}
for any sequence of cubes $(W_n)_{n\in\N}$ with $W_n\rightarrow\R^d$.
Since the Nerve Theorem (see, e.g., p.~76 in \citep{Edelsbrunner}) states that $\beta_p(Z^{W_n})=\beta_p\big(\Delta_{W_n}\big)$ almost surely, the assertion follows.
\end{proof}

The following remark provides a criterion under which condition \eqref{BoolM_BettiInt} is always satisfied for every $p\in\N_0$.

\begin{Bem}
Assume that the Boolean model has bounded grains, meaning that there exists an $R>0$ with
\begin{align*}
    \Theta\big( \big\{ K\in\mathcal{K}^d : K\subseteq B(\0,R) \big\} \big) \,=\, 1,
\end{align*}
where $B(x,r)$ denotes the closed ball with center $x$ and radius $r$.
In this situation, we obtain almost surely
\begin{align*}
\pi(V) \,&=\, \gamma \int_{\R^d} \int_{\mathcal{K}^d} \1\big\{ V\cap (y+K)\neq\emptyset \big\} \; \Theta (\d K) \; \d y \\
&\leq\, \gamma \int_{\R^d} \int_{\mathcal{K}^d} \1\big\{ B(\0,R)\cap B(y,R)\neq\emptyset \big\} \; \Theta (\d K) \; \d y \\
&=\, \gamma (2R)^d\, |B(\0,1)|,
\end{align*}
which means that the mixing parameter $\pi(V)$ is almost surely bounded and therefore, all its moments exist.
\end{Bem}

\section*{Acknowledgements}

The results of this paper stem from the author's PhD thesis \cite{Pabst.Thesis}.
The author is grateful to Daniel Hug, the supervisor of the PhD thesis, for the support and insightful feedback throughout the research.
This work was partially supported by the Deutsche Forschungsgemeinschaft (DFG, German Research Foundation) through the SPP 2265, under grant numbers HU 1874/5-1 and ME 1361/16-1.

\bigskip
\bigskip
\bigskip
\bigskip
\bigskip

\bibliographystyle{abbrv}
\bibliography{RandomConnectionModel_final}

\end{document}